\documentclass[11pt]{article}
\usepackage{epsfig,amssymb}

\begin{document}

\title{
Hexagonal circle patterns with constant intersection angles and
discrete Painlev\'e and Riccati equations}

\author{{\Large Agafonov S.I.\footnote{\tt sagafonov@rusfund.ru} } \& {\Large Bobenko A.I.\footnote{\tt bobenko@math.tu-berlin.de} }\\
\\
Fachbereich Mathematik, Technische Universit\"at Berlin \\ Strasse
des 17. Juni 136, 10623 Berlin, Germany \\ }
\date{}
\maketitle

\unitlength=1mm

\newtheorem{theorem}{Theorem}
\newtheorem{proposition}{Proposition}
\newtheorem{lemma}{Lemma}
\newtheorem{corollary}{Corollary}
\newtheorem{definition}{Definition}

\pagestyle{plain}

\begin{abstract}

\noindent Hexagonal circle patterns with constant intersection angles mimicking holomorphic
 maps $z^{c }$ and ${\rm log}(z)$  are studied.
It is shown that the corresponding
 circle patterns are immersed and described by
special separatrix solutions of  discrete  Painlev\'e and Riccati
equations. The general solution of the Riccati equation is
expressed in terms of the hypergeometric function. Global
properties of these solutions, as well as of the discrete $z^c$
and ${\rm log}(z)$, are established.
\end{abstract}

\section{Introduction. Hexagonal circle patterns and $z^c$}

The theory of circle patterns is a rich fascinating area having
its origin in the classical theory of circle packings. Its fast
development in recent years is caused by the mutual influence and
interplay of ideas and concepts from discrete geometry, complex
analysis and the theory of integrable systems.

The progress in this area  was initiated by Thurston's idea
\cite{T},\cite{MR} of approximating the Riemann mapping by circle
packings. Classical circle packings consisting of disjoint open
disks were later generalized to circle patterns where the disks
may overlap (see for example \cite{H}). Different underlying
combinatorics were considered. Circle patterns with the
combinatorics of the square grid  were introduced in
\cite{Schramm}; hexagonal circle patterns were studied in
\cite{BHS} and \cite{BH}.

The striking analogy between circle patterns and the classical
analytic function theory is underlined by such facts as the
uniformization theorem
 concerning circle packing realizations of
cell complexes with prescribed combinatorics \cite{BS}, a discrete
maximum principle and Schwarz's lemma \cite{R}, rigidity
properties \cite{MR},\cite{H} and a discrete Dirichlet principle
\cite{Schramm}.

The convergence of discrete conformal maps represented by circle
packings was proven in \cite{RS}. For prescribed regular
combinatorics this result was refined. $C^{\infty}$-convergence
for hexagonal packings  is shown in  \cite{HS}. The uniform
convergence for circle patterns with the combinatorics of the
square grid and orthogonal neighbouring circles was established in
\cite{Schramm}.

The approximation issue naturally leads to the question about
analogs to standard holomorphic functions. Computer experiments
give evidence for their existence \cite{DS},\cite{TH}, however not
very much is known.  For circle packings with hexagonal
combinatorics the only explicitly described
 examples are Doyle spirals  \cite{Doy},\cite{BDS}, which are  discrete
analogs of exponential maps, and conformally symmetric packings,
which are analogs of a quotient of Airy functions \cite{BHConf}.
For patterns with overlapping circles more explicit examples are
known: discrete versions of ${\rm exp} (z)$, ${\rm erf}(z)$
\cite{Schramm},  $z^c$, $ {\rm log} (z)$ \cite{AB} are constructed
for patterns with underlying combinatorics of the square grid;
$z^c$, ${\rm log}(z)$ are also described for hexagonal patterns
\cite{BHS}, \cite{BH}.

It turned out that an effective approach to the description  of
circle patterns is given by the theory of integrable systems (see
\cite{BPD},\cite{BHS},\cite{BH}). For example, Schramm's circle
patterns are governed by a difference equation which is the
stationary Hirota equation (see \cite{Schramm}). This approach
proved to be especially useful for the construction of discrete
$z^c$ and ${\rm log}(z)$ in
\cite{AB},\cite{BHS},\cite{BH},\cite{BPD} with the aid of some
isomonodromy problem. Another connection with the theory of
discrete integrable equations was revealed in
\cite{AB},\cite{A},\cite{A1}: embedded circle patterns are
described by special solutions of discrete Painlev\'e II and
discrete Riccati equations.

In this paper we carry the results of \cite{AB} for square grid
combinatorics over to hexagonal circle patterns  with constant
intersection angles introduced in \cite{BH}.

Hexagonal combinatorics are obtained on a sub-lattice of ${\mathbb
Z}^3$ as follows: consider the subset
 $$H=\{(k,l,m)\in {\mathbb Z}^3:\ |k+l+m|\le1 \}$$ and
join by edges those vertices of $H$ whose $(k,l,m)$-labels differ
by 1 only in one component. The obtained quadrilateral lattice
$QL$ has two types of vertices: for $k+l+m=0$ the corresponding
vertices have 6 adjacent edges, while the vertices with $k+l+m=\pm
1$ have only 3. Suppose that the vertices with 6 neighbors
correspond to centers of circles in the complex plane $\mathbb C$
and the vertices with 3 neighbors correspond to intersection
points of circles with the centers in neighboring vertices. Thus
we obtain a circle pattern with hexagonal combinatorics.

\begin{figure}[th]
\begin{center}
\epsfig{file=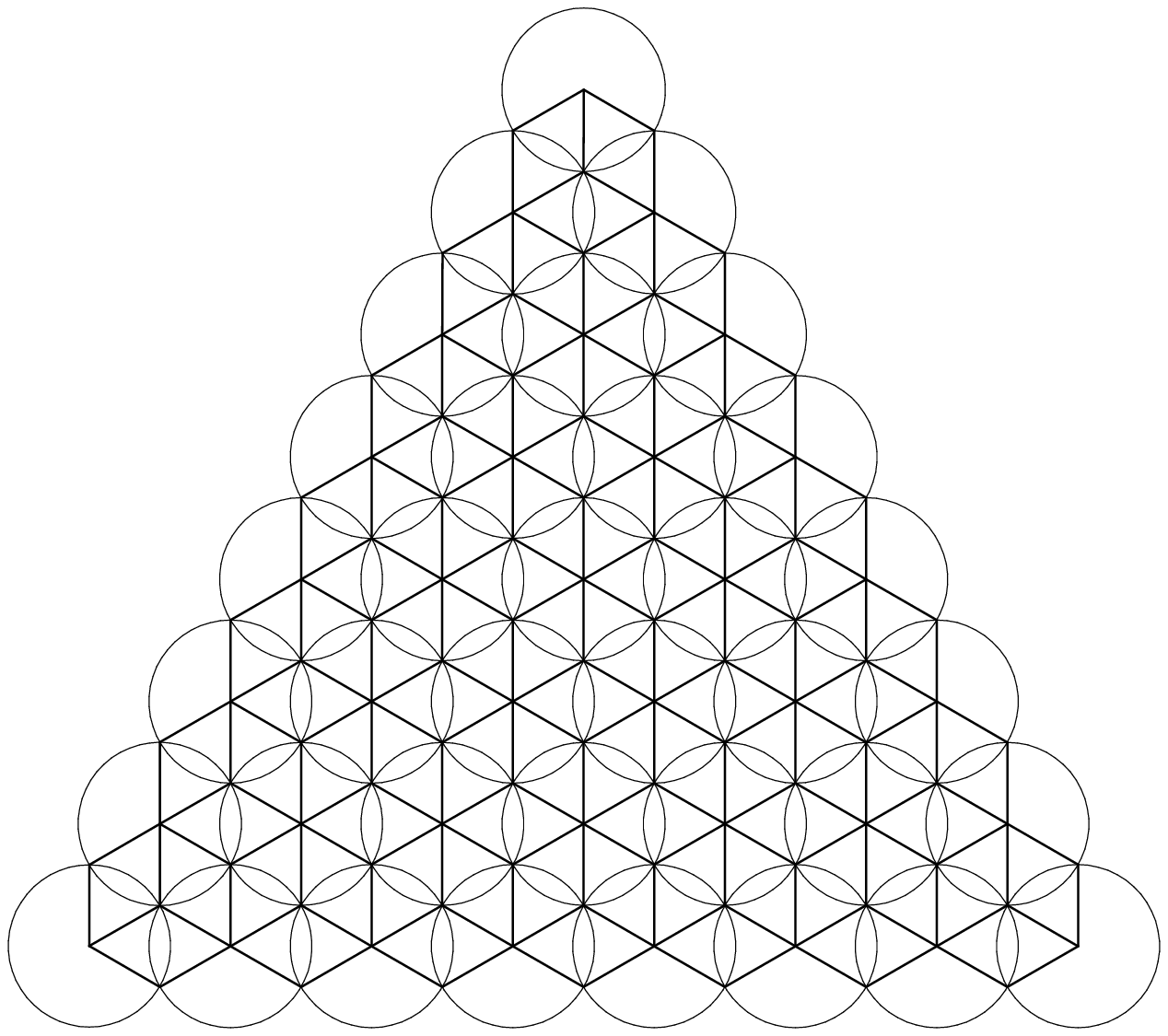,width=45mm}
\begin{picture}(20,50)
\put(2,22){\vector(1,0){15}} \put(7,27){ \it  \huge $Z^{3/2}$}
\end{picture} \ \ \ \ \
\epsfig{file=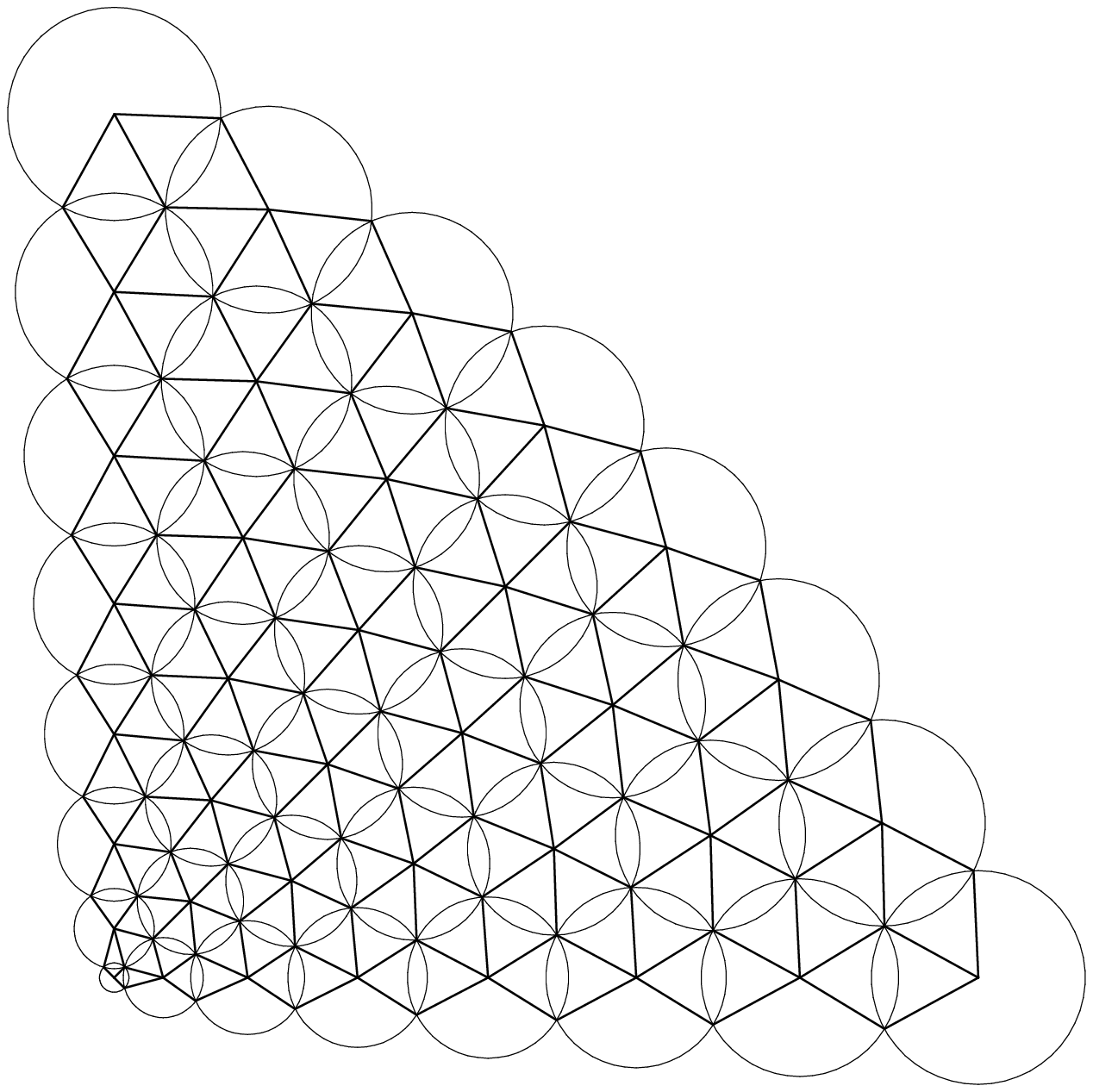,width=45mm}
\caption{Hexagonal circle
patterns as a discrete conformal map.} \label{CPMap}
\end{center}
\end{figure}

Circle patterns where the intersection angles are constant for
each of 3 types of (quadrilateral) faces (see Fig.\ref{CPMap})
were introduced in \cite{BH}. A special case of such circle
patterns mimicking holomorphic map $z^c$ and ${\rm log}(z)$ is
given by the restriction to an $H$-sublattice of a special
isomonodromic solution of some
 {\it integrable system} on the lattice ${\mathbb Z}^3$. Equations
for the field variable $z:{{\mathbb Z}^3}\to {\mathbb C}$ of this
system are:
$$
q(z_{k,l,m},z_{k,l+1,m},z_{k-1,l+1,m},z_{k-1,l,m})=e^{-2i\alpha
_1}, $$
\begin{equation}
\label{skew-cross-ratios}
q(z_{k,l,m},z_{k,l,m-1},z_{k,l+1,m-1},z_{k,l+1,m})=e^{-2i\alpha
_2},
\end{equation}
$$ q(z_{k,l,m},z_{k+1,l,m},z_{k+1,l,m-1},z_{k,l,m-1})=e^{-2i\alpha
_3}, $$ where $\alpha _i > 0$ satisfy $\alpha _1+\alpha
_2+\alpha
 _3=\pi$
and
 $$
q(z_1,z_2,z_3,z_4)=\frac{(z_1-z_2)(z_3-z_4)} {(z_2-z_3)(z_4-z_1)}
$$ is the cross-ratio of elementary quadrilaterals of the image of
${{\mathbb Z}^3}$. Equations (\ref{skew-cross-ratios}) mean
 that the cross-ratios of images of
faces of elementary cubes are constant for each type of faces,
while the restriction $\alpha _1+\alpha_2+\alpha _3=\pi$ ensures
their compatibility.

The isomonodromic problem for this system (see section
\ref{secMon} for the details, where we present the necessary
results from  \cite{BH} ) specifies the non-autonomous constraint
$$ c
z_{k,l,m}=2k\frac{(z_{k+1,l,m}-z_{k,l,m})(z_{k,l,m}-z_{k-1,l,m})}{z_{k+1,l,m}-z_{k-1,l,m}}+
$$
\begin{equation}\label{constraint}
\qquad \qquad
2l\frac{(z_{k,l+1,m}-z_{k,l,m})(z_{k,l,m}-z_{k,l-1,m})}{z_{k,l+1,m}-z_{k,l-1,m}}+
\end{equation}
$$ \qquad \qquad
2m\frac{(z_{k,l,m+1}-z_{k,l,m})(z_{k,l,m}-z_{k,l,m-1})}{z_{k,l,m+1}-z_{k,l,m-1}},
$$ which is compatible with (\ref{skew-cross-ratios})  (this
constraint in the two-dimensional case with $c=1$ first appeared
in \cite{NC}). In particular, a solution to
(\ref{skew-cross-ratios}),(\ref{constraint}) in the subset
\begin{equation}\label{domain}
Q=\{(k,l,m)\in {{\mathbb Z}^3}| \ k\ge0,\ l\ge0, \ m\le 0  \}
\end{equation}
is uniquely determined by its values
$$
 z_{1,0,0},\ z_{0,1,0}, \
z_{0,0,-1}.
$$
  Indeed, the  constraint (\ref{constraint}) gives
$z_{0,0,0}=0$ and  defines $z$ along the coordinate axis
$(n,0,0)$, $(0,n,0)$, $(0,0,-n)$. Then all other $z_{k,l,m}$ with
$(k,l,m)\in Q$ are calculated through the cross-ratios
(\ref{skew-cross-ratios}).

\begin{proposition}\cite{BH}\label{circularity}
The solution $z:Q \to {\mathbb C}$ of the system
(\ref{skew-cross-ratios}),(\ref{constraint}) with the initial data
\begin{equation}\label{general_initial}
  z_{1,0,0}=1,\ z_{0,1,0}=e^{i\phi}, \
z_{0,0,-1}=e^{i\psi}
\end{equation}
determines a circle pattern. For all $(k,l,m)\in Q$ with even
$k+l+m$ the points $z_{k\pm 1,l,m}$, $z_{k,l\pm 1,m}$, $z_{k,l,m\pm 1}$
lie on a circle with the center $z_{k,l,m}$, i.e. all elementary
quadrilaterals of the $Q$-image are of kite form.
\end{proposition}

 Moreover, equations (\ref{skew-cross-ratios}) (see lemma \ref{geom-kite}
 in section \ref{secEuc}) ensure that
for the points $z_{k,l,m}$ with $k+l+m=\pm 1$ where 3 circles meet
intersection angles are $\alpha_i$ or $\pi-\alpha _i$, $i=1,2,3$
 (see Fig.\ref{CPMap} where the isotropic case
 $\alpha_i=\pi/3$  of regular and $Z^{3/2}$-pattern are
 shown).

According to proposition \ref{circularity}, the discrete map
$z_{k,l,m}$, restricted on $H$, defines a circle pattern with
circle centers $z_{k,l,m}$ for $k+l+m=0$, each circle intersecting
6 neighboring circles. At each intersection points three circles
meet.

However, for most initial data $\phi, \psi \in {\mathbb R}$, the
behavior of the obtained circle pattern is quite irregular: inner
parts of different elementary quadrilaterals intersect and circles
overlap. Define $Q_H=Q\cap H$.

\begin{definition}
\label{def} \cite{BH} The hexagonal circle pattern $Z^c$,
  $0<c<2$  with  intersection
 angles $\alpha_1,\alpha_2,\alpha _3$, $\alpha _i > 0$,  $\alpha_1+\alpha_2+\alpha_3=\pi$
 is the  solution $z:Q \to {\mathbb C}$ of
(\ref{skew-cross-ratios})  subject to (\ref{constraint}) and with
the initial data
 \begin{equation}\label{initial}
  z_{1,0,0}=1,\
  z_{0,1,0}=e^{ic  (\alpha _2+\alpha _3)  }, \
z_{0,0,-1}=e^{ic \alpha _3}
\end{equation}
restricted to $Q_H$.
\end{definition}

\begin{definition}
A discrete  map $z:Q_H \to {\bf C}$ is called an immersion, if
inner parts of adjacent elementary quadrilaterals are disjoint.
\end{definition}
The main result of this paper is the following theorem.
\begin{theorem}\label{main}
The hexagonal  $Z^{c}$ with constant positive intersection angles
and $0<c <2$ is an immersion.
\end{theorem}

 The proof of this
property follows from an analysis of the geometrical properties of
the corresponding circle patterns and analytical properties of the
corresponding discrete Painlev\'e and Riccati equations.

The crucial step is to consider equations for the radii of the
studied circle patterns in the whole $Q$-sublattice  with even
$k+l+m$. In section \ref{secEuc}, these equations are derived and
the geometrical property of immersedness is reformulated as the
positivity of the solution to these equations. Using discrete
Painlev\'e and Riccati equations in section \ref{secImm} we
present the proof of the existence of a positive solution and thus
complete the proof of immersedness. In section \ref{secGen}, we
discuss possible generalizations and corollaries of the obtained
results. In particular, circle patterns $Z^2$ and $\rm Log$ with
both square grid and hexagonal combinatorics are considered. It is
also proved that they are immersions.

\section{Discrete $Z^c$ via a monodromy problem}\label{secMon}

Equations (\ref{skew-cross-ratios})  have the Lax representation
\cite{BH}: $$ \Phi
_{k+1,l,m}(\mu)=L^{(1)}(e,\mu)\Phi_{k,l,m}(\mu), $$
\begin{equation}\label{Lax}
\Phi _{k,l+1,m}(\mu)=L^{(2)}(e,\mu)\Phi_{k,l,m}(\mu),
\end{equation}
$$\Phi _{k,l,m+1}(\mu)=L^{(3)}(e,\mu)\Phi_{k,l,m}(\mu),$$ where
$\mu$ is the spectral parameter and $\Phi(\mu):{{\mathbb Z}^3}\to
{\rm GL}(2,{\mathbb C})$ is the wave function. The matrices
$L^{(n)}$ are defined on the edges $e=(p_{out},p_{in})$ of
${{\mathbb Z}^3}$ connecting two neighboring vertices and oriented
in the direction of increasing $k+l+m$:
\begin{equation}\label{matrix}
L^{(n)}(e,\mu) =\left(
\begin{array}{cc}
1 & z_{in}-z_{out}\\ \mu\frac{\Delta _n}{z_{in}-z_{out}}& 1
\end{array}
\right),
\end{equation}
with parameters $\Delta_n$ fixed for each type of edges. The
zero-curvature condition on the faces of elementary cubes of
${{\mathbb Z}^3}$ is equivalent to equations
(\ref{skew-cross-ratios}) with $\Delta _n =e^{i\delta _n}$ for
properly chosen $\delta _n$. Indeed, each elementary quadrilateral
of ${{\mathbb Z}^3}$ has two consecutive positively oriented pairs
of edges $e_1,e_2$ and $e_3,e_4$. Then the compatibility condition
$$ L^{(n_1)}(e_2) L^{(n_2)}(e_1)= L^{(n_2)}(e_4) L^{(n_1)}(e_3)$$
is exactly one of the equations (\ref{skew-cross-ratios}). This
Lax representation is a generalization of the one found in
\cite{NC} for the square lattice.

A solution $z:{{\mathbb Z}^3}\to {\mathbb C}$ of equations
(\ref{skew-cross-ratios}) ia called {\it isomonodromic} if there
exists a wave function $\Phi(\mu):{\rm \bf Z^3}\to {\rm GL}(2,{\bf
C})$ satisfying (\ref{Lax}) and the following linear differential
equation in $\mu$:
\begin{equation}\label{monodromy}
\frac{d}{d \mu}\Phi _{k,l,m}(\mu)=A_{k,l,m}(\mu)\Phi
_{k,l,m}(\mu),
\end{equation}
where $A_{k,l,m}(\mu)$ are some $2\times 2$ matrices meromorphic
in $\mu$, with the order and position of their poles being
independent of $k,l,m$.

The simplest non-trivial isomonodromic solutions satisfy the
constraint:
$$bz^2_{k,l,m}+ c
z_{k,l,m}+d=2(k-a_1)\frac{(z_{k+1,l,m}-z_{k,l,m})(z_{k,l,m}-z_{k-1,l,m})}{z_{k+1,l,m}-z_{k-1,l,m}}+
$$
\begin{equation}\label{constraintMob}
\qquad \qquad
2(l-a_2)\frac{(z_{k,l+1,m}-z_{k,l,m})(z_{k,l,m}-z_{k,l-1,m})}{z_{k,l+1,m}-z_{k,l-1,m}}+
\end{equation}
$$ \qquad \qquad
2(m-a_3)\frac{(z_{k,l,m+1}-z_{k,l,m})(z_{k,l,m}-z_{k,l,m-1})}{z_{k,l,m+1}-z_{k,l,m-1}}.
$$

\begin{theorem}\cite{BH}
Let $z:{\bf Z^3}\to {\bf C}$ be an  isomonodromic solution  to
(\ref{skew-cross-ratios}) with the matrix $A_{k,l,m}$ in
(\ref{monodromy}) of the form
\begin{equation}\label{matrix-form}
A_{k,l,m}(\mu)=\frac{C_{k,l,m}}{\mu}+\sum_{n=1}^3\frac{B^{(n)}_{k,l,m}}{\mu-\frac{1}{\Delta_n}}
\end{equation}
with $\mu$-independent matrices $C_{k,l,m}$, $B^{(n)}_{k,l,m}$ and
normalized by ${\rm tr}\ A_{0,0,0}(\mu)=0$. Then these matrices
have the following form: $$ C_{k,l,m} =\frac{1}{2}\left(
\begin{array}{cc}
-bz_{k,l,m}-c/2 & bz^2_{k,l,m}+cz_{k,l,m}+d\\
 b&  bz_{k,l,m}+c/2
\end{array}
\right) $$
$$ B^{(1)}_{k,l,m}
=\frac{k-a_1}{z_{k+1,l,m}-z_{k-1,l,m}}\left(
\begin{array}{cc}
z_{k+1,l,m}-z_{k,l,m} &
(z_{k+1,l,m}-z_{k,l,m})(z_{k,l,m}-z_{k-1,l,m})\\
 1& z_{k,l,m}-z_{k-1,l,m}
\end{array}
\right) +\frac{a_1}{2}I$$ $$ B^{(2)}_{k,l,m}
=\frac{l-a_2}{z_{k,l+1,m}-z_{k,l-1,m}}\left(
\begin{array}{cc}
z_{k,l+1,m}-z_{k,l,m} &
(z_{k,l+1,m}-z_{k,l,m})(z_{k,l,m}-z_{k,l-1,m})\\
 1& z_{k,l,m}-z_{k,l-1,m}
\end{array}
\right) +\frac{a_2}{2}I$$ $$ B^{(3)}_{k,l,m}
=\frac{m-a_3}{z_{k,l,m+1}-z_{k,l,m-1}}\left(
\begin{array}{cc}
z_{k,l,m+1}-z_{k,l,m} &
(z_{k,l,m+1}-z_{k,l,m})(z_{k,l,m}-z_{k,l,m-1})\\
 1& z_{k,l,m}-z_{k,l,m-1}
\end{array}
\right) +\frac{a_3}{2}I$$ and $z_{k,l,m}$ satisfies
(\ref{constraintMob}).

Conversely, any solution $z:{{\mathbb Z}^3}\to {\mathbb C}$ to the
system (\ref{skew-cross-ratios}),(\ref{constraintMob}) is
isomonodromic with $A_{k,l,m}(\mu)$ given by the formulas above.
\end{theorem}
The special case $b=a_1=a_2=a_3=0$ with shift $z\to z-d/c$ implies
(\ref{constraint}).

\section{Euclidean description of hexagonal circle patterns}
\label{secEuc}

In this section we describe the circle pattern $z^c$ in terms of
the radii of the circles. Such characterization proved to be quite
useful for the circle patterns with combinatorics of the square
grid \cite{AB},\cite{A}. In what follows, we say that the triangle
$(z_1,z_2,z_3)$ has {\it positive (negative) orientation} if
$$\frac{z_3-z_1}{z_2-z_1}=\left|\frac{z_3-z_1}{z_2-z_1}\right|e^{i\phi}
\ \ {\rm with} \ 0\le  \phi \le \pi \  \ \ (-\pi <\phi < 0).$$

\begin{lemma}\label{geom-kite}
Let $q(z_1,z_2,z_3,z_4)=e^{-2i\alpha }$, $0<\alpha <\pi$.
\begin{itemize}
\item If $|z_1-z_2|=|z_1-z_4|$ and the triangle $(z_1,z_2,z_4)$ has
positive orientation then $|z_3-z_2|=|z_3-z_4|$ and the angle
between $[z_1,z_2]$ and $[z_2,z_3]$ is $(\pi-\alpha )$.
\item If $|z_1-z_2|=|z_1-z_4|$ and the triangle $(z_1,z_2,z_4)$ has
negative orientation then $|z_3-z_2|=|z_3-z_4|$ and the angle
between $[z_1,z_2]$ and $[z_2,z_3]$ is $\alpha$.
\item If the angle
between $[z_1,z_2]$ and $[z_1,z_4]$ is $\alpha$ and the triangle
$(z_1,z_2,z_4)$ has positive orientation then
$|z_3-z_2|=|z_1-z_2|$ and   $|z_3-z_4|=|z_4-z_1|$.
\item If the angle
between $[z_1,z_2]$ and $[z_1,z_4]$ is $(\pi- \alpha)$ and the
triangle $(z_1,z_2,z_4)$ has negative orientation then
$|z_3-z_2|=|z_1-z_2|$ and   $|z_3-z_4|=|z_4-z_1|$.
\end{itemize}
\end{lemma}
Lemma \ref{geom-kite} and Proposition \ref{circularity} imply that
each elementary quadrilateral of the studied circle pattern  has
one of the forms enumerated in the lemma.

Proposition \ref{circularity} allows us to introduce the radius
function
\begin{equation} \label{r-def}
r({\scriptscriptstyle K,L,M})=|z_{k,l,m}-z_{k \pm 1,l,m}|=|z_{k,l
\pm 1,m}-z_{k,l,m}|=|z_{k,l,m}-z_{k,l,m\pm 1}|,
\end{equation}
where $(k,l,m)$ belongs to the sublattice of $Q$ with even $k+l+m$ and
$(K,L,M)$ label this sublattice:
\begin{equation} \label{sublattice}
 K=k-\frac{k+l+m}{2}, \ \
L=l-\frac{k+l+m}{2}, \ \  M=m-\frac{k+l+m}{2}.
\end{equation}
The function $r$ is defined on the sublattice
$$
\tilde Q=\{(K,L,M)\in {\bf Z^3}| L+M\le 0,\ \ M+K\le 0,\ \
K+L\ge 0 ) \}
$$
corresponding to $Q$. Consider this function on
$$
\tilde Q_H=\{(K,L,M)\in {\bf Z^3}| K\ge 0,\ \ L\ge 0,\ \ M\le 0,
\ \ K+L+M=0,+1\}.
$$

\begin{figure}[th]
\begin{center}
\epsfig{file=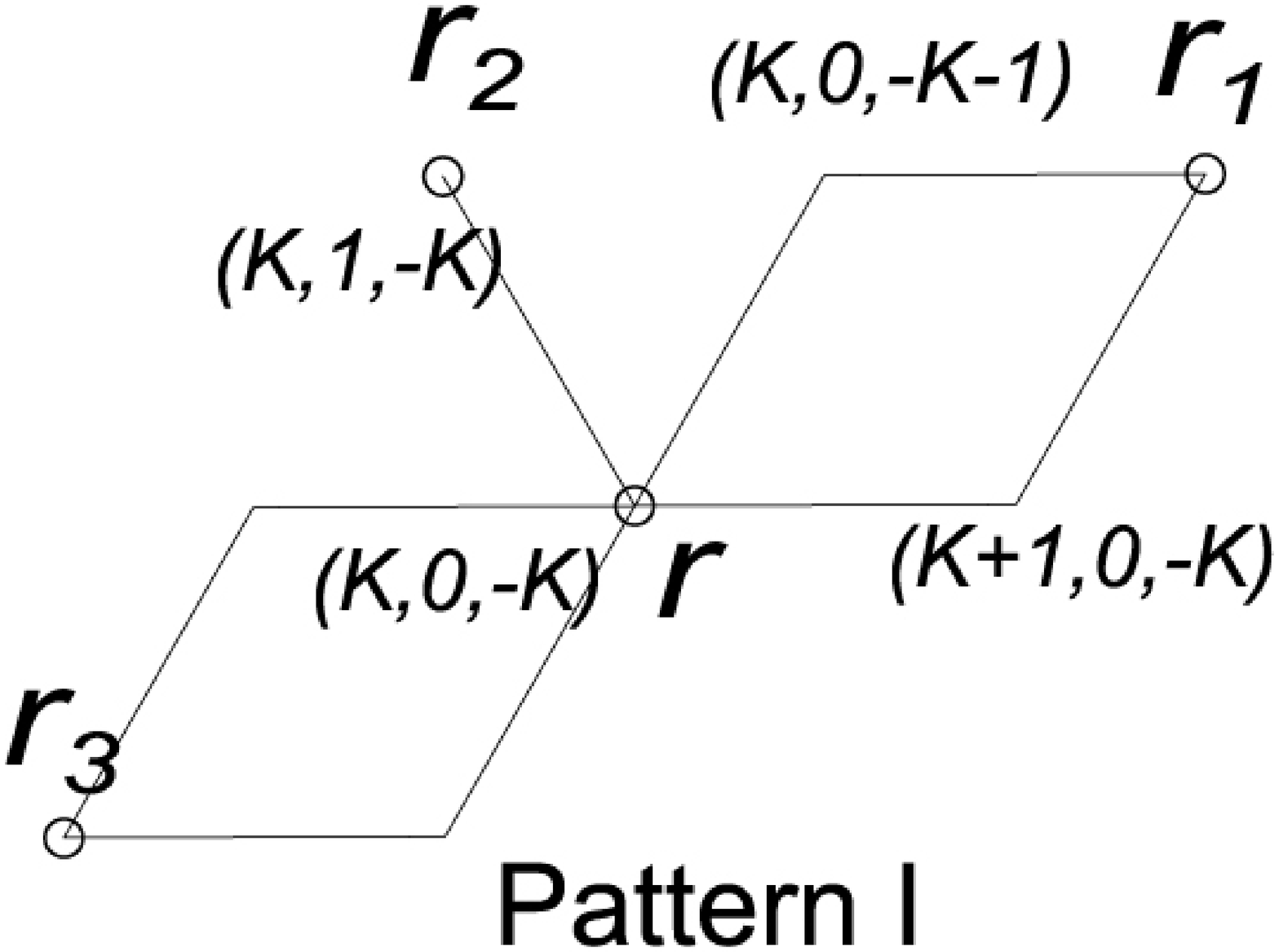,height=30mm}\ \ \
\epsfig{file=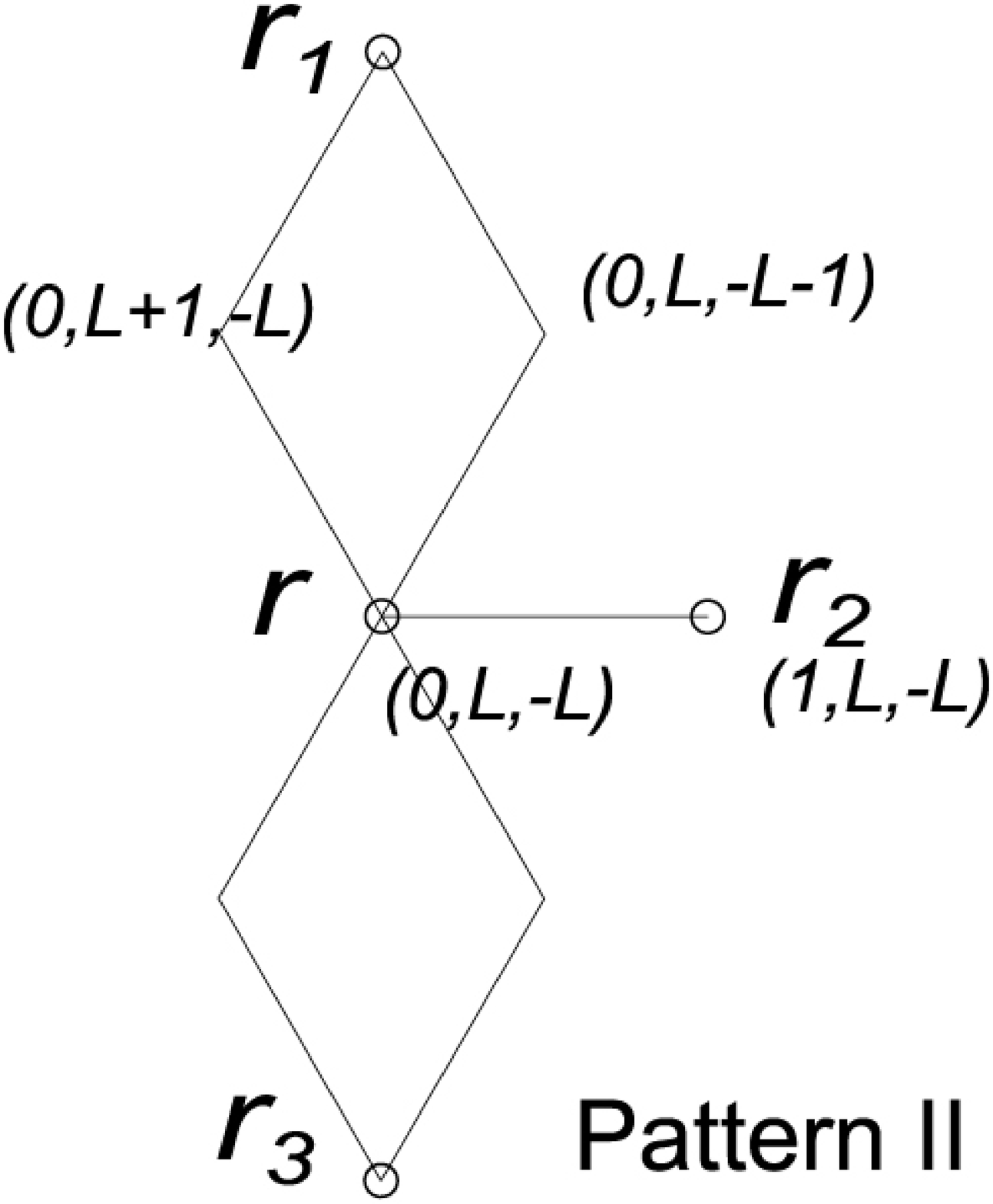,height=42mm} \ \ \
\epsfig{file=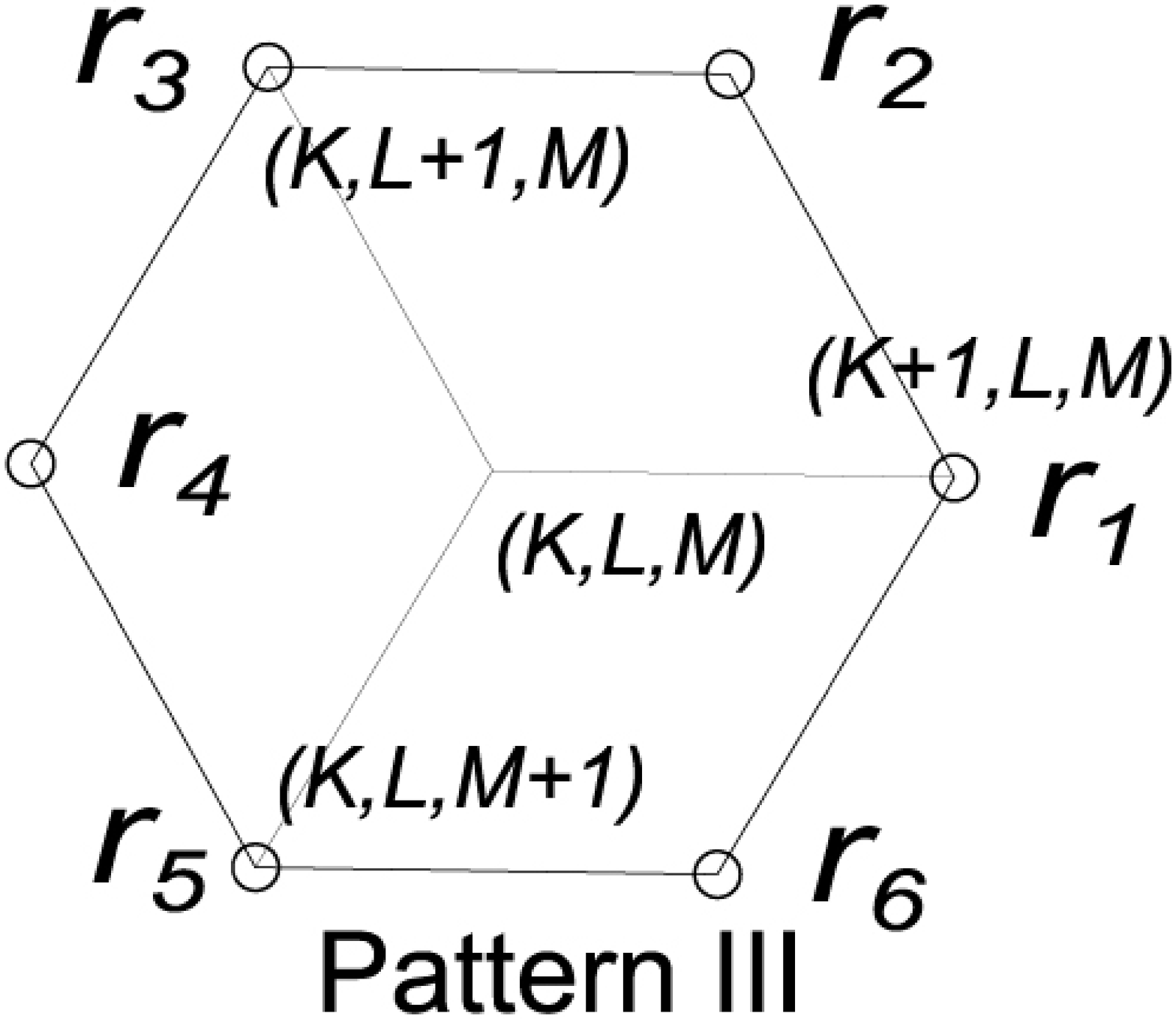,height=32mm} \ \ \ \ \
\epsfig{file=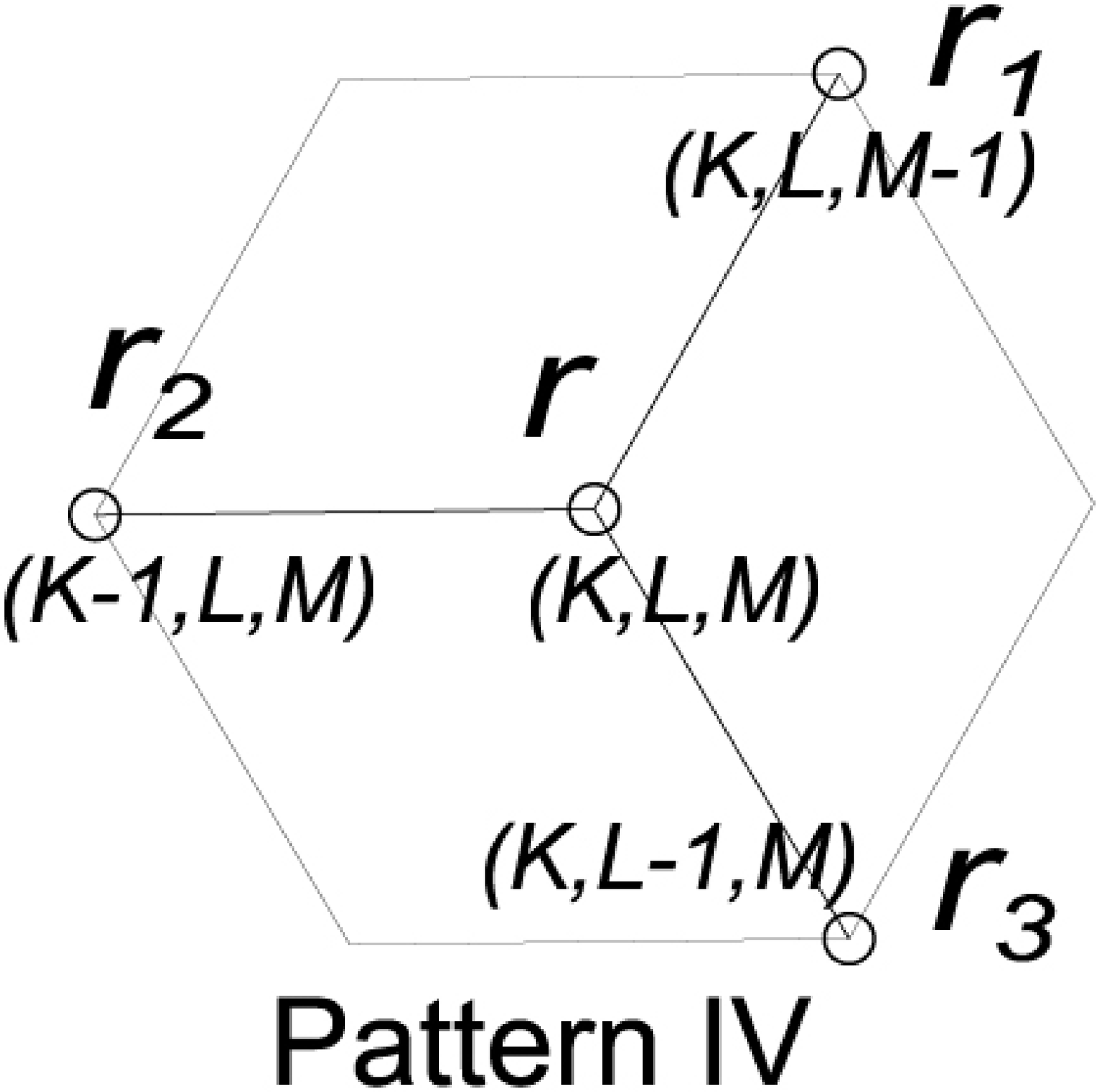,width=32mm}
 \caption{Equation patterns.}
\label{Rb}
\end{center}
\end{figure}

\begin{theorem} \label{rimmersion}
Let the solution $z:Q_H \to {\mathbb C}$ of the system
(\ref{skew-cross-ratios}),(\ref{constraint}) with initial data
(\ref{general_initial}) be an immersion.  Then function
$r({\scriptscriptstyle K,L,M}): \tilde Q_H \to {{\mathbb R}_+}$,
defined by (\ref{r-def}), satisfies the following equations:
\begin{equation}\label{Rkm}
(r_1+r_2)(r^2-r_2r_3+r(r_3-r_2)\cos \alpha_i )+
(r_3+r_2)(r^2-r_2r_1+ r(r_1-r_2)\cos \alpha_i )=0
\end{equation}
on the patterns of  type I and II as in Fig.\ref{Rb}, with $i=3$ and $i=2$ respectively,
\begin{equation}\label{HexEq}
{\scriptstyle (L+M+1)}\frac{r_4-r_1}{r_4+r_1}+{\scriptstyle
(M+K+1)}\frac{r_6-r_3}{r_6+r_3}+{\scriptstyle
(K+L+1)}\frac{r_2-r_5}{r_2+r_5}=c -1
\end{equation}
on the patterns of type III and
\begin{equation}\label{TriEq}
r(r_1\sin \alpha_3+r_2 \sin \alpha_1 +r_3\sin \alpha_2 )=
r_1r_2\sin \alpha_2+r_2r_3\sin \alpha_3+r_3r_1\sin \alpha_1
\end{equation}
on the patterns of type IV. Conversely, $r({\scriptscriptstyle
K,L,M}): \tilde Q_H \to {{\mathbb R}_+}$ satisfying equations
(\ref{Rkm}),(\ref{HexEq}),(\ref{TriEq}) is the radius function an
immersed hexagonal circle pattern with constant intersection
angles (i.e. corresponding to some immersed solution $z:Q_H \to
{\mathbb C}$ of (\ref{skew-cross-ratios}),(\ref{constraint})),
which is determined by $r$ uniquely.
\end{theorem}
\noindent {\it Proof:} The map $z_{k,l,m}$ is an immersion if and only
if all triangles $(z_{k,l,m},z_{k+1,l,m},z_{k,l,m-1})$,\\
$(z_{k,l,m},z_{k,l,m-1},z_{k,l+1,m})$ and
$(z_{k,l,m},z_{k+1,l,m},z_{k,l+1,m})$  of elementary
quadrilaterals of the map $z_{k,l,m}$ have the same orientation
(for brevity we call it the orientation of the quadrilaterals).

\smallskip

\begin{figure}[th]
\begin{center}
\epsfig{file=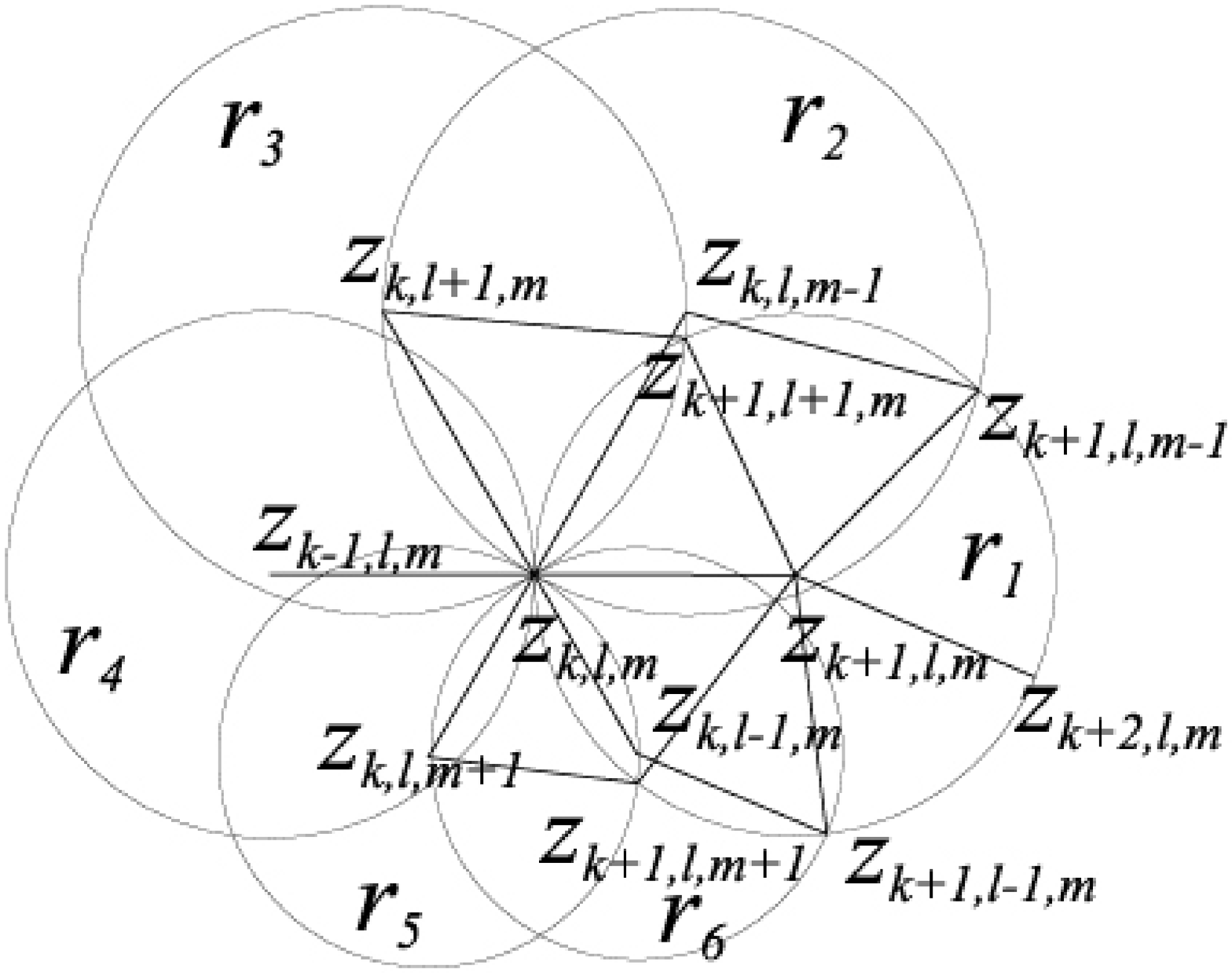,width=65mm}\ \
\epsfig{file=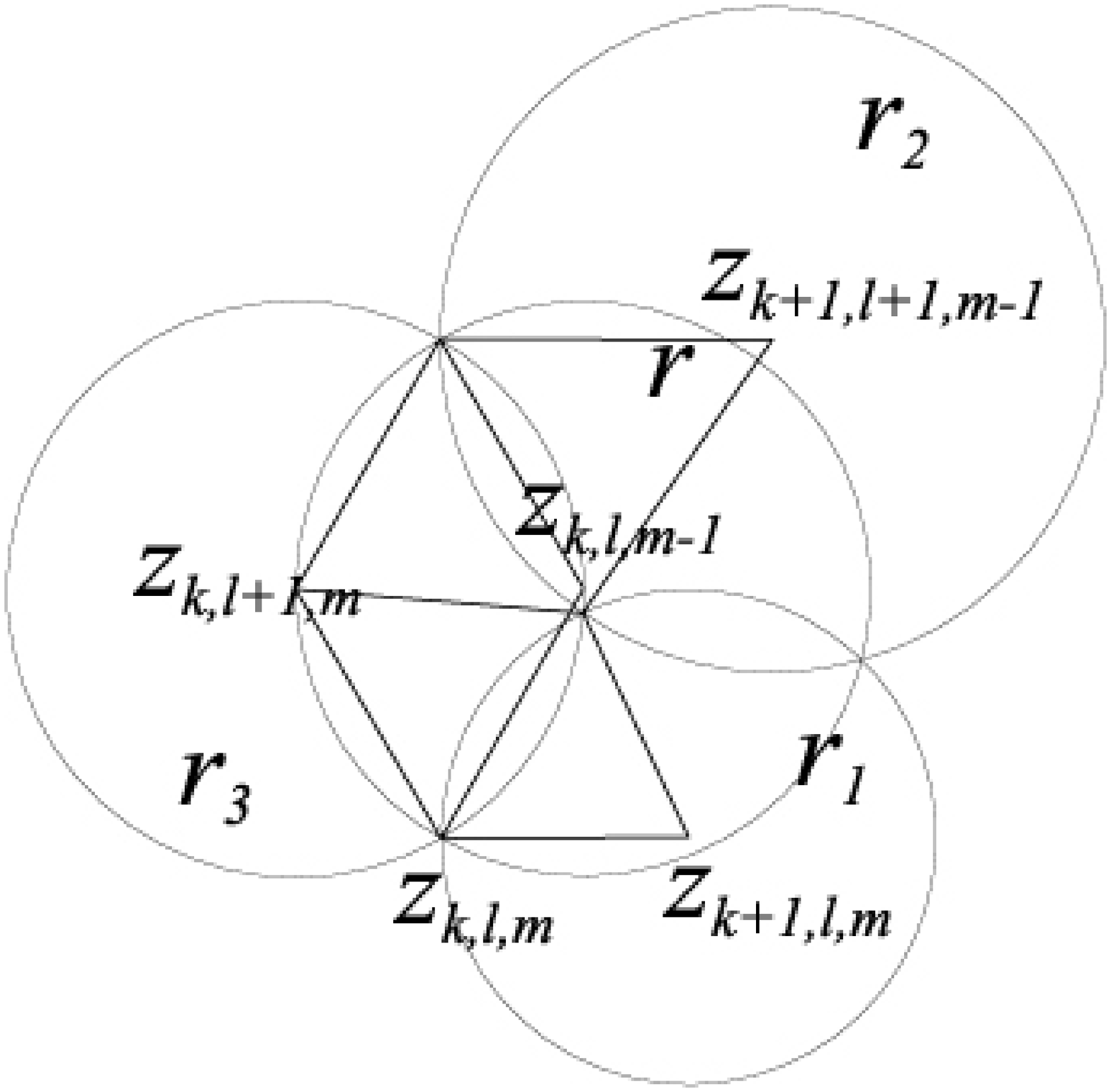,width=55mm}
 \caption{Circles.}
\label{SuperCross}
\end{center}
\end{figure}

\noindent {\it Necessity:} To get equation (\ref{HexEq}), consider
the configuration of two star-like figures with centers at
$z_{k,l,m}$ with $k+l+m=1 \ ({\rm mod} \ 2)$ and at $z_{k+1,l,m}$,
connected by five edges in the $k$-direction as shown on the left
part of Fig.\ref{SuperCross}.
 Let $r_i,\  i=1,...,6$  be the radii of the
circles with the centers at the vertices neighboring  $z_{k,l,m}$
as in Fig.\ref{SuperCross}. As follows from Lemma \ref{geom-kite},
the vertices $z_{k,l,m}$, $z_{k+1,l,m}$ and $z_{k-1,l,m}$ are
collinear. For immersed $z^c$, the vertex $z_{k,l,m}$ lies between
$z_{k+1,l,m}$ and $z_{k-1,l,m}$. Similar facts are true also for
the $l$- and $m$-directions. Moreover, the orientations of
elementary quadrilaterals with the vertex $z_{k,l,m}$ coincides
with one of the standard lattice. Lemma \ref{geom-kite} defines
all angles at $z_{k,l,m}$ of these quadrilaterals. Equation
(\ref{constraint}) at $(k,l,m)$  gives $z_{k,l,m}$:
 $$ z_{k,l,m}=
\frac{2e^{is}}{c
}\left(k\frac{r_1r_4}{r_1+r_4}+l\frac{r_3r_6}{r_3+r_6}e^{i(\alpha_2+\alpha_3)}+m\frac{r_2r_5}{r_2+r_5}e^{i(\alpha_1+\alpha_2
+2 \alpha_3)}\right),$$ where
$e^{is}=(z_{k+1,l,m}-z_{k,l,m})/r_1$. Lemma \ref{geom-kite} allows
one to compute $z_{k+1,l,m-1}$, $z_{k+1,l+1,m}$, $z_{k+1,l,m+1}$
and $z_{k+1,l-1,m}$ using the form of quadrilaterals (they are
shown in Fig. \ref{SuperCross}). Now equation (\ref{constraint})
at $(k+1,l,m)$  defines $z_{k+2,l,m}$. Condition
$|z_{k+2,l,m}-z_{k+1,l,m}|=r_1$ with the labels (\ref{sublattice})
yields equation (\ref{HexEq}).

For $l=0$ values $z_{k+1,0,m}$, $z_{k+2,0,m}$, $z_{k+1,0,m-1}$ and
the equation for the cross-ratio with $\alpha_3$  give the radius
$R$ with the center at $z_{k+2,0,m-1}$. Note that for $l=0$ the
term with $r_6$ and $r_5$ drops out of equation  (\ref{HexEq}).
Using this equation and the permutation $R\to r_1$, $r_1\to r$,
$r_2\to r_2$, $r_5 \to r_3$, one gets equation (\ref{Rkm}) with
$i=3$. The equation for pattern II is derived similarly.

To derive (\ref{TriEq}), consider the figure on the right part of
Fig.\ref{SuperCross} where $k+l+m=1 \ ({\rm mod} \ 2)$ and $r_1$,
$r_2$, $r_3$ and $r$ are the radii of the circles with the centers
at $z_{k+1,l,m}$, $z_{k+1,l+1,m-1}$, $z_{k,l+1,m}$ and
$z_{k,l,m-1}$, respectively. Elementary geometrical considerations
and Lemma \ref{geom-kite} applied to the forms of the shown
quadrilaterals give equation (\ref{TriEq}).

\medskip

\noindent {\bf Remark.} Equation (\ref{TriEq}) is derived for
$r=r({\scriptscriptstyle K,L,M})$, $r_1=r({\scriptscriptstyle
K,L,M-1})$, $r_2=r({\scriptscriptstyle K-1,L,M})$,
$r_3=r({\scriptscriptstyle K,L-1,M+1})$. However  it  holds true
also for
 $r_1=r({\scriptscriptstyle K,L,M+1})$,
$r_2=r({\scriptscriptstyle K+1,L,M})$, $r_3=r({\scriptscriptstyle
K,L+1,M+1})$ since it gives the radius of the circle through the
three intersection points of the circles with radii $r_1$, $r_2$,
$r_3$ intersecting at prescribed angles as shown in the right part
of Fig.\ref{SuperCross}. Later, we refer to this equation also for
this pattern.

\medskip

\begin{figure}[th]
\begin{center}
\epsfig{file=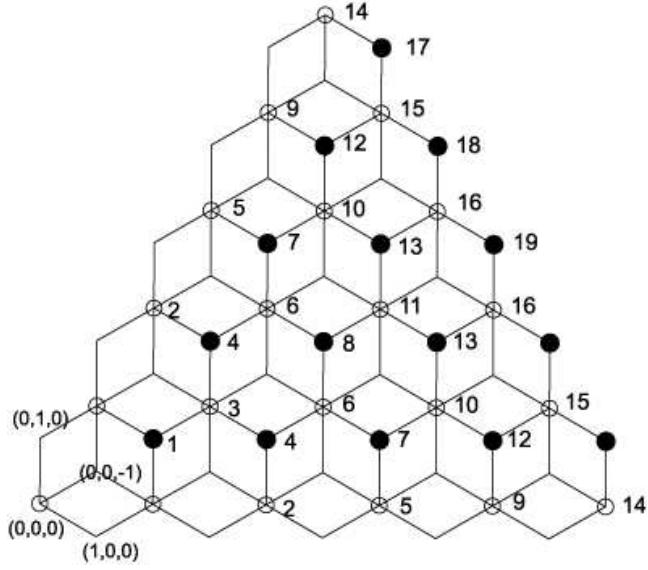,width=85mm}
 \caption{Computing $r$ in $\tilde Q_H$.}
\label{Lattice}
\end{center}
\end{figure}

\noindent {\it Sufficiency:} Now let $r({\scriptscriptstyle
K,L,M}): \tilde Q_H \to {\bf R_+}$ be some positive solution to
(\ref{Rkm}),(\ref{HexEq}),(\ref{TriEq}). We can re-scale it so
that $r({0,0,0})=1$. Starting with $r(1,0,-1)$ and $r(0,1,-1)$ one
can compute $r$ everywhere in $\tilde Q_H$: $r$ in a "black"
vertex (see Fig.\ref{Lattice}) is computed from (\ref{HexEq}).
(Note that only $r$ at "circled" vertices is used: so to compute
$r_{1,1,-1}$ one needs only $r(1,0,-1)$ and $r(0,1,-1)$.) The
function $r$ in "white" vertices on the border $\partial \tilde
Q_H=\{(K,0,-K)|K\in {\bf N}\}\cup \{(0,L,-L)|L\in {\bf N}\}$ is
given by (\ref{Rkm}). Finally, $r$ in "white" vertices in
$Q^{int}_H=Q_H \setminus \partial \tilde Q_H$ is computed from
(\ref{TriEq}). In  Fig.\ref{Lattice} labels show the order of
computing $r$.
\begin{lemma}\label{extra}
Any solution $r({\scriptscriptstyle
K,L,M}): \tilde Q_H \to {\bf R}$  to
 (\ref{Rkm}),(\ref{HexEq}),(\ref{TriEq}) with $0\le c \le 2$,
 which is positive for inner
 vertices of $\tilde Q_H$   defines some $z_{k,l,m}$
 satisfying (\ref{skew-cross-ratios}) in $Q$. Moreover,
 all the triangles $(z_{k,l,m},z_{k+1,l,m},z_{k,l,m-1})$,
$(z_{k,l,m},z_{k,l,m-1},z_{k,l+1,m})$ and
$(z_{k,l,m},z_{k+1,l,m},z_{k,l+1,m})$ have positive orientation.
\end{lemma}
{\it Proof of the lemma:} One can place the circles with radii
$r({\scriptscriptstyle K,L,M})$ into the complex plane $\mathbb C$
in the way prescribed by the hexagonal combinatorics and the
intersection angles. Taking the circle centers and the
intersection points of neighboring circles, one recovers
$z_{k,l,m}$ for $k+l+m=0,\pm 1 $ up to a translation and rotation.
Reversing the arguments used in the derivation of
(\ref{Rkm}),(\ref{HexEq}),(\ref{TriEq}), one observes from the
forms of the quadrilaterals that equations
(\ref{skew-cross-ratios}) are satisfied. Now using
(\ref{skew-cross-ratios}), one recovers $z$ in the whole $Q$.
Equation (\ref{TriEq}) ensures that the radii $r$ remain positive,
which implies the positive orientations of the triangles
$(z_{k,l,m},z_{k+1,l,m},z_{k,l,m-1})$,
$(z_{k,l,m},z_{k,l,m-1},z_{k,l+1,m})$ and
$(z_{k,l,m},z_{k+1,l,m},z_{k,l+1,m})$.

\medskip

Consider a solution $z:Q \to {\mathbb C}$ of the system
(\ref{skew-cross-ratios}),(\ref{constraint}) with initial data
(\ref{general_initial}), where $\phi$ and $\psi$ are chosen so
that the triangles $(z_{0,0,0)},z_{1,0,0},z_{0,0,-1})$ and
$(z_{0,0,0)},z_{0,0,-1},z_{0,1,0})$ have positive orientations and
satisfy  conditions $r(1,0,-1)=|z_{1,0,-1}-z_{1,0,0}|$ and
$r(0,1,-1)=|z_{0,1,-1}-z_{0,0,-1}|$. The map $z_{k,l,m}$ defines
circle pattern due to Proposition \ref{circularity} and  coincides
with the map defined by Lemma \ref{extra} due to the uniqueness of
the solution uniqueness. Q.E.D.

\medskip

\noindent Since the cross-ratio equations and the constraint are
compatible, the equations for the radii are also compatible.
Starting with $r(0,0,0)$, $r(1,0,-1)$ and $r(0,1,-1)$, one can
compute $r({\scriptscriptstyle K,L,M})$ everywhere in $\tilde Q$.

\begin{lemma} \label{planes}
Let a solution $r({\scriptscriptstyle K,L,M}):\tilde Q \to
{\mathbb R} $ of (\ref{Rkm}),(\ref{HexEq}),(\ref{TriEq}) be
positive in the planes
 given by equations $K+M=0$ and $L+M=0$ then it is positive everywhere
 in $\tilde Q$.
\end{lemma}
\noindent {\it Proof:} As follows from equation (\ref{TriEq}), $r$
is positive for positive $r_i$, $i=1,2,3.$ As $r$ at $(K,K,-K)$,
$(K+1,K,-K-1)$ and $(K,K+1,-K-1)$ is positive, $r$ at  $(K,K,K-1)$
is also positive. Now starting from $r$ at $(K,K,-K-1)$ and having
$r>0$ at $(N,K+1,-K-1)$ and $(N,K,-K)$, one obtains positive $r$
at $(N,K,-K-1)$ for $0\le N <K$ by the same reason. Similarly, $r$
at $(K,N,-K-1)$ is positive. Thus  from positive $r$ at the planes
$K+M=0$ and $L+M=0$, we get positive $r$ at the planes $K+M=-1$
and $L+M=-1$. Induction completes the proof.

\medskip

\begin{lemma} \label{lines}
Let a solution $r({\scriptscriptstyle K,L,M}):\tilde Q \to
{\mathbb R} $ of (\ref{Rkm}),(\ref{HexEq}),(\ref{TriEq}) be
positive in the lines
 parameterized by $n$ as  $(n,0,-n)$ and $(0,n,-n)$. Then it is positive
 in the border   planes of $\tilde Q$ specified by  $K+M=0$ and $L+M=0$.
\end{lemma}
\noindent {\it Proof:} We prove this lemma for $K+M=0$. For the
other border plane it is proved similarly. Equation (\ref{HexEq})
for $(K,L,-K-1)$ gives
\begin{equation} \label{square}
r_2=r_5\frac{(2L+c)r_1+(2K+c)r_4}{(2K+2-c)r_1+(2L+2-c)r_4}
\end{equation}
therefore $r_2$ is positive provided $r_1$, $r_5$ and $r_4$ are
positive. For $K=L$ it reads as
 \begin{equation} \label{lin}
r_2=r_5\frac{(2K+c)}{(2K+2-c)}.
\end{equation}
It allows us to compute recursively $r$ at $(K,K,-K)$ starting
with $r$ at $(0,0,0)$. Obviously, $r>0$ for $(K,K,-K)$ if $r>0$ at
$(0,0,0)$. This property together with the condition $r>0$ at
$(n,0,-n)$ imply the conclusion of the lemma since equation
(\ref{square}) gives $r$ everywhere in the border   plane of
$\tilde Q$ specified by  $K+M=0$.

\medskip

\noindent Lemmas \ref{planes} and \ref{lines} imply  that the
circle pattern $z^c$ is an immersion if $r>0$ at $(N,0,-N)$ and
$(0,N,-N)$.

\section{Proof of the main theorem. Discrete Painlev\'e and Riccati
 equations} \label{secImm}

In this section, we prove that all $r(n,0,-n), \ \forall n\in
{\mathbb N}$ are positive only for the initial data $z_{1,0,0}=1$,
$z_{0,0,-1}=e^{c\alpha_3}$. For the line $r(0,n,-n)$ the proof is
the same. Our strategy is as follows: first, we prove the
existence of an initial value $z_{0,0,-1}$ such that $r(n,0,-n)>0,
\ \forall n\in {\mathbb N}$, finally we will show that this value
is unique and is $z_{0,0,-1}=e^{c\alpha_3}$.

\begin{proposition} \label{P}
Suppose the equation
\begin{equation}
\label{dPII} (n+1)(x_n^2-1)\left(\frac{x_{n+1}+{x_n}/{\varepsilon}
}{\varepsilon+x_nx_{n+1}}\right)-
 n(1-{x_n^2}/{\varepsilon^2})\left(\frac{x_{n-1}+\varepsilon x_n}{\varepsilon+x_{n-1}x_n}\right) =
c x_n \frac{\varepsilon ^2 -1}{2\varepsilon^2},
\end{equation}
where $\varepsilon =e^{i\alpha_3}$, has a unitary solution
$x_n=e^{i\beta_n}$ in the sector $0<\beta_n < \alpha _3$. Then
$r(n,0,-n)$, $n\ge 0$ is positive.
\end{proposition}
\noindent {\it Proof:} For $z_{1,0,0}=1$ and unitary $z_{1,0,-1}$,
the equation for the cross-ratio with $\alpha _3$ and
(\ref{constraint}) reduce to (\ref{dPII}) with unitary
$x_n^2=(z_{n,0,-n-1}-z_{n,0,n})/(z_{n+1,0,-n}-z_{n,0,n})$. Note
that for $n=0$ the term with $x_{-1}$ drops out of (\ref{dPII});
therefore the solution for $n>0$ is determined by $x_0$ only. The
condition $0<\beta_n < \alpha _3$ means that all triangles
$(z_{n,0,-n},z_{n+1,0,n},z_{n,0,-n-1})$ have positive orientation.
Hence $r(n,0,-n)$ are all positive. Q.E.D.

\medskip

\noindent {\bf Remark.} Equation (\ref{dPII}) is a special {\it
discrete Painlev\'e equation}. For a more general reduction of
cross-ratio equation see \cite{N}. The case $\varepsilon=i$,
corresponding to the orthogonally intersecting circles, was
studied in detail in \cite{AB}. Here we generalize these results
to the case of arbitrary unitary $\varepsilon$. Below we omit the
index of $\alpha$ so that $\varepsilon =e^{i\alpha}$.

\begin{theorem} \label{exist}
There exists a unitary solution  $x_n=e^{i\beta_n}$ to (\ref{dPII})
in the sector $0<\beta_n < \alpha $.
\end{theorem}

\noindent {\it Proof:} Equation (\ref{dPII}) allows us to
represent $x_{n+1}$ as a function of $n,x_{n-1}$ and $x_n$:
$x_{n+1}=\Phi (n,x_{n-1},x_n)$. $\Phi(n,u,v)$ maps the torus
$T^2=S^1\times S^1=\{(u,v)\in {\mathbb C}:\ |u|=|v|=1\}$ into
$S^1$ and has the following properties:
\begin{itemize}
\item For all $n\in {\mathbb N}$ it is a continuous map on $A_I\times \bar A_I$
where $A_I=\{e^{i\beta}: \beta \in (0,\alpha)\}$ and $\bar A_I$ is the
      closure of $A_I$ . Values of $\Phi$ on
 the border of $A_I\times \bar A_I$ are defined by continuity:
      $\Phi(n,u,\varepsilon)=-1$, $\Phi(n,u,1)=-\varepsilon$.
\item For $(u,v)\in A_I\times A_I$ one has $\Phi(n,u,v)\in A_I\cup A_{II} \cup
A_{IV}$, where $A_{II}=\{e^{i\beta}: \beta \in (\alpha,\pi]\}$ and
$A_{IV}=\{e^{i\beta}: \beta \in [\alpha-\pi,0)\}$. I.e., $x$
cannot jump in one step from $A_I$ into $A_{III}=\{e^{i\beta}:
\beta \in (-\pi,\alpha-\pi )\}$.
\end{itemize}
Let $x_0=e^{i\beta_0}$. Then $x_n=x_n(\beta_0)$. Define
$S_n=\{\beta_0: \ x_k(\beta_0)\in \bar A_I \ \forall \ 0\le k \le
n\}$. Then $S_n$ is a closed set since $\Phi$ is continuous on
$A_I\times \bar A_I$. As a closed subset of a segment it is a
collection of disjoint segments $S_n^l$.
\begin{lemma} \label{segments}
There exists sequence $\{S_n^{l(n)}\}$ such that:
\begin{itemize}
\item  $S_n^{l(n)}$  is mapped by $x_n(\beta_0)$ onto
$\bar A_I$,
\item  $S_{n+1}^{l(n+1)}\subset  S_n^{l(n)}$.
\end{itemize}
\end{lemma}
The lemma is proved by induction. For $n=0$ it is trivial. Suppose
it holds for $n$. As $S_n^{l(n)}$  is mapped by $x_n(\beta_0)$
onto $\bar A_I$, continuity considerations and
$\Phi(n,u,\varepsilon)=-1$, $\Phi(n,u,1)=-\varepsilon$ imply:
$x_{n+1}(\beta _0)$ maps $S_n^{l(n)}$ onto $ A_I\cup A_{II} \cup
A_{IV}$ and at least one of the segments $S_{n+1}^l\subset
S_n^{l(n)}$ is mapped into $\bar A_I$. This proves the lemma.

\medskip

Since the segments of $\{S_n^{l(n)}\}$ constructed in lemma
\ref{segments} are nonempty, there exists $\bar \beta_0\in S_n$
for all $n\ge 0$. For this $\bar \beta _0$, the value  $x_n(\bar
\beta _0)$ is not on the border of $\bar A _0$ since then
$x_{n+1}(\beta_0)$ would jump out of $\bar A_I$. Q.E.D..

\medskip

\begin{figure}[th]
\begin{center}
\epsfig{file=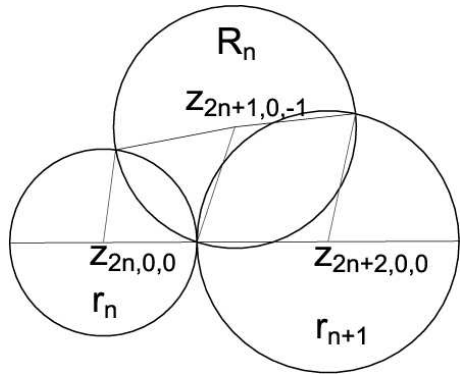,width=60mm} \caption{Circles on the border.
 } \label{RicC}
\end{center}
\end{figure}

Let $r_n$ and $R_{n}$ be the radii of the circles of the circle
patterns defined by $z_{k,l,m}$ with the centers at $z_{2n,0,0}$
and $z_{2n+1,0,-1}$ respectively. Constraint (\ref{constraint})
gives
$$r_{n+1}=\frac{2n+c}{2(n+1)-c}r_{n}$$ which is exactly formula
(\ref{lin}). From elementary geometric considerations (see Fig.
\ref{RicC})
one gets
$$
R_{n+1}=\frac{r_{n+1}-R_n \cos \alpha}{R_n-r_{n+1}\cos
\alpha}r_{n+1}
$$
(recall that $\alpha =\alpha _3$).
 Define $$ p_n=\frac{R_{n}}{r_n}, \ \ \
g_n(c)=\frac{2n+c}{2(n+1)-c}$$ and denote $t=\cos \alpha$. Now,
the equation for the radii $R,\ r$ takes the form:
\begin{equation}\label{Riccati}
p_{n+1}=\frac{g_n(c)-tp_n}{p_n-tg_n(c)}.
\end{equation}

\medskip

\noindent {\bf Remark.} Equation (\ref{Riccati}) can be seen as a
discrete version of a {\it  Riccati} equation. This is motivated
by the following properties:
\begin{itemize}
\item the cross-ratio of each four-tuple of its solutions is constant
as $p_{n+1}$ is a M\"obius transform of $p_n$,
\item the general solution is expressed in terms of solutions of
some linear equation: the standard Ansatz
\begin{equation}\label{ansatz}
p_n=\frac{y_{n+1}}{y_n}+tg_n(c )
\end{equation}
transforms (\ref{Riccati}) into
\begin{equation}\label{linear}
y_{n+2}+t(g_{n+1}(c)+1)y_{n+1}+(t^2-1)g_n(c )y_n=0.
\end{equation}

\end{itemize}

\medskip

\noindent As follows from Theorem \ref{exist}, Proposition \ref{P}
and Lemma \ref{lines}, equation (\ref{Riccati}) has a positive
solution. One may conjecture that there is only one initial value
$p_0$ such that $p_n>0,\ \forall n\in {\mathbb N}$ from the
consideration of the asymptotics. Indeed, $g_n(c)\to 1$ as $n\to
\infty$, and the general solution of equation (\ref{linear}) with
limit values of coefficients is $y_n=c_1(-1)^n(1+t)^n+c_2(1-t)^n$.
Thus $p_n=\frac{y_{n+1}}{y_n}+tg_n(c)\to -1$ for $c_1\ne 0$.
However $c_1,c_2$ define only the asymptotics of a solution. To
relate it to the initial value $p_0$ is a more difficult problem.
Fortunately, it is possible to find the general solution to
(\ref{linear}).

\begin{proposition}\label{solution}
The general solution to (\ref{linear}) is
\begin{equation}\label{series}
y_n=\frac{\Gamma(n+\frac{1}{2})}{\Gamma (n+1-\frac{c
}{2})}(c_1\lambda_1 ^{n+1-c /2}F(\frac{3-c
}{2},\frac{c -1}{2},\frac{1}{2}-n ,z_1)+
\end{equation}
$$ +c_2\lambda_2 ^{n+1-c/2}F(\frac{3-c
}{2},\frac{c -1}{2},\frac{1}{2}-n ,z_2) $$ where $\lambda _1=-t-1,
\ \lambda _2=1-t, \ z_1 =(t-1)/2,\ z_2= -(1+t)/2$ and $F$ is the
hypergeometric function.
\end{proposition}
{\it Proof:} The solution was found by a slightly modified {\it
symbolic method} (see \cite{Boole} for the method description and
\cite{A1} for the detail). Here, $F(a,b,c,z)$ denotes the standard
hypergeometric function which is a solution of the hypergeometric
equation
\begin{equation}\label{Gauss-eq}
z(1-z)F_{zz}+[c-(a+b+1)z]F_z-a b F=0
\end{equation}
holomorphic at $z=0$.  Due to linearity, the general solution of
(\ref{linear}) is given by a superposition of any two linearly
independent solutions. Direct computation with the series
representation of the hypergeometric function
\begin{equation} \label{hyper}
F(1-\frac{\gamma-1}{2},\frac{\gamma-1}{2},
1-\left(x+\frac{\gamma-1}{2}\right),z)
=1+z\frac{(1-\frac{\gamma-1}{2})(\frac{\gamma-1}{2})}{(1-(x+\frac{\gamma-1}{2}))}
+\ldots+ \end{equation}
$$
z^k\frac{ \left[(1-\frac{\gamma-1}{2})(2-\frac{\gamma-1}{2})
\ldots(k-\frac{\gamma-1}{2})\right]\left[(\frac{\gamma-1}{2})
(1+\frac{\gamma-1}{2})\ldots
(k-1\frac{\gamma-1}{2})\right]}{(1-(x\frac{\gamma-1}{2}))
\ldots(k-(x+\frac{\gamma-1}{2}))}+\ldots
$$
 shows that each summand  in (\ref{series})  satisfies equation
(\ref{linear}). To finish the proof of Proposition \ref{solution},
one has to show that the particular solutions with $c_1=0,\ c_2\ne
0$ and $c_1 \ne 0,\ c_2=0$
 are linearly independent. This fact follows from
\begin{lemma}\label{as-lin}
As $n\to \infty $, function (\ref{series})  has the asymtotics
\begin{equation}\label{eq-as-lin}
y_n \simeq (n+1-\gamma /2)^{\frac{\gamma -1}{2}}(c_1\lambda_1
^{n+1-\gamma /2}+c_2\lambda_2 ^{n+1-\gamma /2}).
\end{equation}
 \end{lemma}
For $n\to \infty$ the series representation (\ref{hyper}) implies
$F(\frac{3-\gamma }{2},\frac{\gamma -1}{2},\frac{1}{2}-n
,z_1)\simeq1$. Stirling's formula
\begin{equation}\label{Stirling}
\Gamma (x) \simeq \sqrt{2\pi }e^{-x}x^{x-\frac{1}{2}}
\end{equation}
yields  the asymptotics of the factor
$\frac{\Gamma(n+\frac{1}{2})}{\Gamma (n+1-\frac{\gamma }{2})}$.
This completes the proof of the lemma and of Proposition
(\ref{solution}).

\begin{proposition}\label{R-positive}
A solution of the discrete Riccati equation (\ref{Riccati}) with
$\alpha \ne \pi/2$  is positive for all $n\ge 0$ if and only if
\begin{equation}\label{p-initial}
p_0=\frac{\sin  \frac{c \alpha}{2}}{\sin \frac{(2-c) \alpha}{2}}.
\end{equation}
\end{proposition}

\noindent {\it Proof:} For positive $p_n$, it is necessary that
$c_1=0$: this follows from asymptotics (\ref{eq-as-lin})
substituted into (\ref{ansatz}). Let us define
\begin{equation}\label{s-def}
s(z)=1+z\frac{(1-\frac{\gamma -1}{2})(\frac{\gamma
-1}{2})}{\frac{1}{2}}+\ldots+z^k\frac{(k-\frac{\gamma
-1}{2})\ldots(1-\frac{\gamma -1}{2})(\frac{\gamma
-1}{2})(k-1+\frac{\gamma
-1}{2})}{k!(k-\frac{1}{2})\ldots\frac{1}{2}}\ldots
\end{equation}
This is the hypergeometric function $F(\frac{3-\gamma
}{2},\frac{\gamma -1}{2},\frac{1}{2}-n ,z)$ with $n=0$. A
straightforward computation with series shows that
\begin{equation}\label{p-0}
p_0=1+\frac{2(\gamma -1)}{2-\gamma }z+\frac{4z(z-1)}{2-\gamma
}\frac{s^\prime (z)}{s(z)}
\end{equation}
where $z=\frac{1+t}{2}$. Note that $p_0$ as a function of $z$
satisfies an ordinary differential equation of first order since
$\frac{s^\prime (z)}{s(z)}$ satisfies the Riccati equation
obtained by a reduction of (\ref{Gauss-eq}). A computation shows
that  $\frac{\sin \frac{\gamma \alpha}{2}}{\sin \frac{(2-\gamma )
\alpha}{2}}$ satisfies the same ordinary differential eqiation.
Since both expression (\ref{p-0}) and (\ref{p-initial}) are equal
to 1 for $z=0$, they coincide everywhere. Q.E.D.

\medskip

\noindent {\it Proof of Theorem \ref{main}:} Proposition
(\ref{R-positive}) implies that the initial $x_0$ for which
(\ref{dPII}) gives positive $r$ is unique and implies the initial
values (\ref{initial}) for $z^c$ if $\alpha _i\ne\pi/2$. For the
case $\alpha =\pi/2$, any solution for (\ref{Riccati}) with
$p_0>0$ is positive. Nevertheless, as was proved in \cite{AB},
$x_0$ is in this case also unique and is specified by
(\ref{p-initial}). Thus for all $n\in {\mathbb N}$ we have
$r(n,0,-n)>0$, $r(0,n,-n)>0$ for the circle pattern $z^c$. Lemmas
\ref{lines} and \ref{planes} complete the proof.

\section{Hexagonal circle patterns $z^2$ and $\rm Log$}
For
$c=2$, formula (\ref{lin}) gives infinite $r(1,1,-1)$. The way
around this difficulty is renormalization $z\to (2-c)z/c$ and a
limit procedure $c\to 2-0$, which leads to the re-normalization of
initial data (see \cite{BH}). As follows from (\ref{p-initial}),
this renormalization implies:
\begin{equation} \label{2-initial}
r(0,0,0)=0, \ r(1,0,-1)=\frac{\sin \alpha _3}{\alpha _3}, \
r(0,1,-1)=\frac{\sin \alpha _2}{\alpha _2}, \ r(1,1,-1)=1.
\end{equation}

\begin{proposition}\label{2-positive}
The solution to (\ref{Rkm}),(\ref{HexEq}),(\ref{TriEq}) with $c=2$
and initial data (\ref{2-initial}) is positive.
\end{proposition}
{\it Proof:} This follows from Lemmas \ref{planes} and \ref{lines}
since Theorem \ref{exist} is true also for the case $c=2$. Indeed,
solution $x_n$ is a continuous function of $c$. Therefore it has a
limit value as $c \to 2-0$ and it lies in the sector $A_I$.

\noindent Lemma \ref{extra} implies that there exists a hexagonal
circle pattern with radius function $r$.
\begin{definition}
The hexagonal circle pattern $Z^2$ has a radius function specified
by Proposition \ref{2-positive}.
\end{definition}
Equations (\ref{Rkm}),(\ref{HexEq}),(\ref{TriEq}) have the
symmetry
\begin{equation} \label{duality}
r \to \frac{1}{r}, \ c\to 2-c,
\end{equation}
which is the {\it duality transformation} (see \cite{BPD}). The
smooth analog $f\to f^*$ for holomorphic functions $f(w),f^*(w)$
is:
$$
\frac{df(w)}{dw}\frac{df^*(w)}{dw}=1.
$$
Note that ${\rm log }^* (w)=w^2/2$.  The hexagonal circle pattern
$\rm Log$ is defined \cite{BH} as a circle pattern dual to $Z^2$.
Discrete $z^2$ and $\rm Log$ are the first two images in Fig.
\ref{CPMapSG}.
\begin{theorem}
The hexagonal circle patterns $Z^2$ and $\rm Log$ are immersions.
\end{theorem}
{\it Proof:} For $z^2$ this follows from Proposition
\ref{2-positive}. Hence the values of $1/r$, where $r$ is radius
function for $z^2$, are positive except for $r(0,0,0)=\infty $.
Lemma \ref{extra} completes the proof.

\section{Concluding remarks} \label{secGen}
In this section we discuss corollaries of the obtained results and
possible generalizations.

\subsection{Square grid circle patterns $z^c$ and $\rm Log$}

\begin{figure}[th]
\begin{center}

\epsfig{file=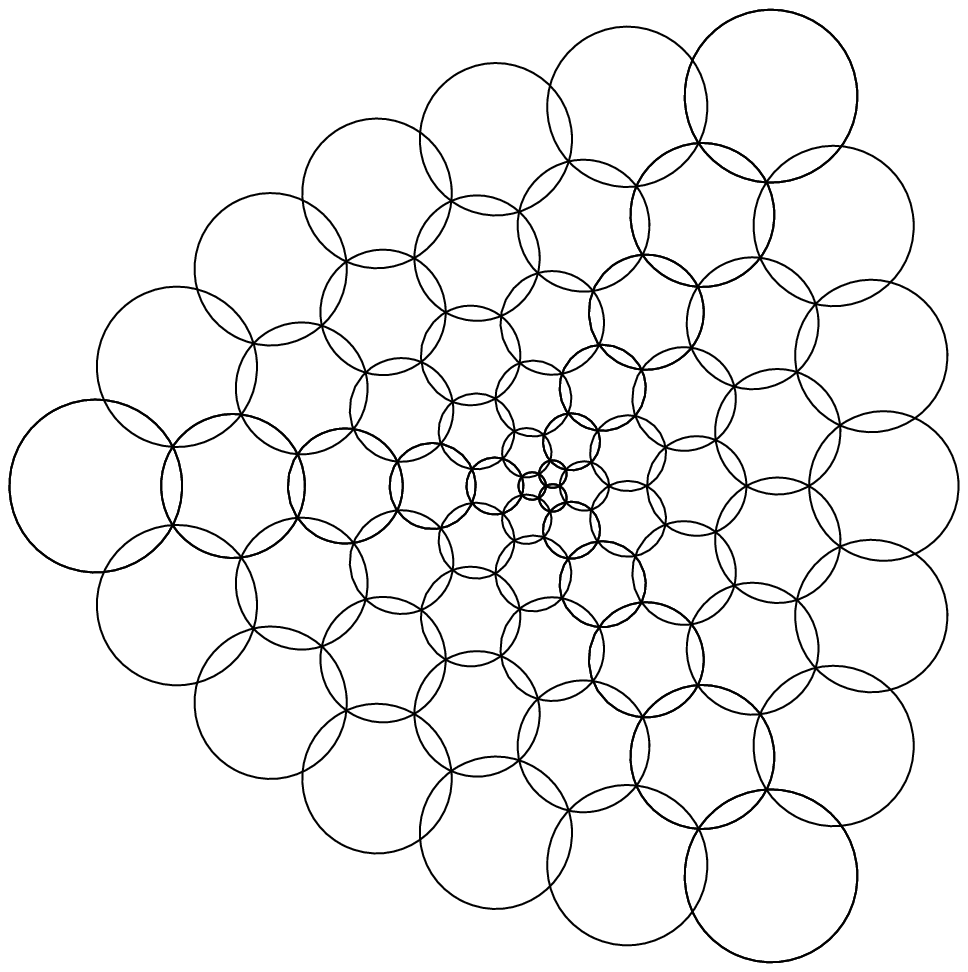,width=45mm}
\epsfig{file=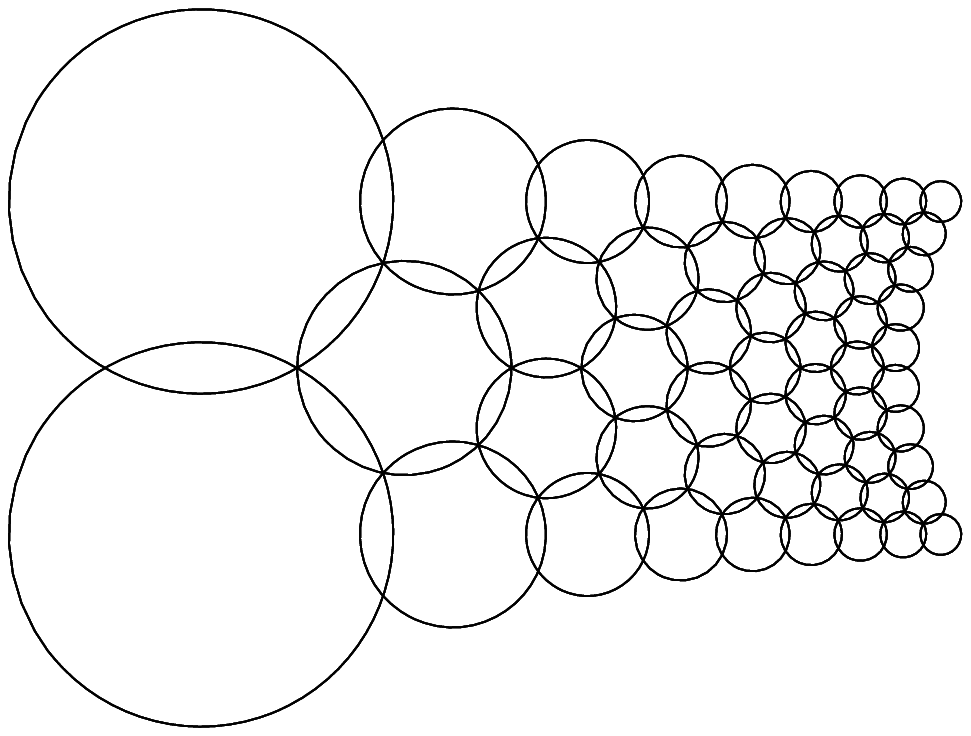,width=50mm}
\epsfig{file=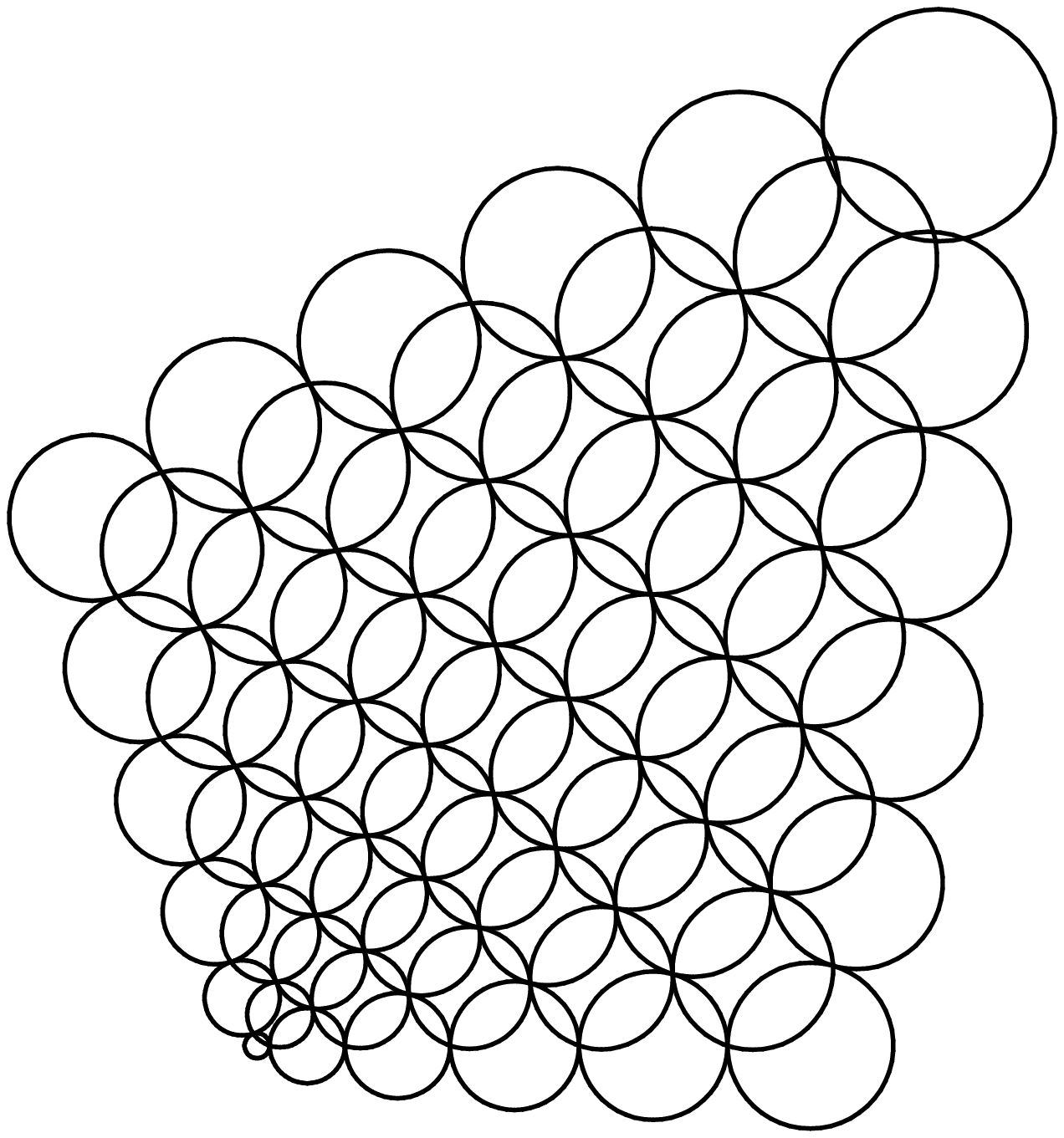,width=40mm}
 \caption{Hexagonal $z^2$, $\rm Log$ and square grid $z^{3/2}$.}
\label{CPMapSG}
\end{center}
\end{figure}

Equations (\ref{skew-cross-ratios}) extend $z_{k,l,m}$
corresponding to the hexagonal $z^2$ and $\rm Log $ from $Q_H$
into the three-dimensional lattice $Q$.  The $r$-function of this
extension satisfies equation (\ref{TriEq}). Consider $z_{k,l,m}$
for the hexagonal $z^c$ and $\rm Log$  restricted to one of the
coordinate planes, e.g. $l=0$. Then Proposition \ref{circularity}
states that $z_{k,0,m}$ defines some circle pattern with
combinatorics of the square grid: each circle has four neighboring
circles intersecting it at angles $\alpha_3$ and $\pi -\alpha _3$.
It is natural to call it {\it square grid} $z^c$ (see the third
image in Fig. \ref{CPMapSG}). Such circle patterns are natural
generalization of those with orthogonal neighboring and tangent
half-neighboring circles introduced and studied in \cite{Schramm}.

\begin{theorem}
Square grid $z^c$, $0< c\le 2$ and   $\rm Log$ are immersions.
\end{theorem}
Proof easily follows from lemma \ref{extra}.

It is interesting to note that the square grid circle pattern
$z^c$ can be obtained from hexagonal one by limit procedure
$\alpha_ 3\to +0$ and by $\alpha _1 \to\pi-\alpha_2$. These limit
circle patterns still can be defined by (\ref{skew-cross-ratios}),
(\ref{constraint}) by imposing the self-similarity condition
$z_{k,l,m}=f_{l,k-m}$.

\subsection{Square grid circle patterns $\rm Erf$}

For square grid combinatorics and $\alpha =\pi/2$, Schramm
\cite{Schramm} constructed circle patterns mimicking holomorphic
function ${\rm  erf} (z)=(2/\pi)\int e^{-z^2}dz$ by giving the
radius function explicitly. Namely, let $n,m$ label the circle
centers so that the pairs of circles $c(n,m)$, $c(n+1,m)$ and
$c(n,m)$, $c(n,m+1)$ are orthogonal and the pairs $c(n,m)$,
$c(n+1,m+1)$ and $c(n,m+1)$, $c(n+1,m)$ are tangent. Then
\begin{equation}\label{erf-r}
r(n,m)=e^{nm}
\end{equation}
 satisfies the equation for a radius function:
\begin{equation} \label{ortho-sg}
R^2(r_1+r_2+r_3+r_4)-(r_2r_3r_4+r_1r_3r_4+r_1r_2r_4+r_1r_2r_3)=0,
\end{equation}
where $R=r(n,m)$, $r_1=r(n+1,m)$, $r_2=r(n,m+1)$, $r_3=r(n-1,m)$,
$r_4=r(n,m-1)$. For square grid circle patterns with intersection
angles $\alpha$ for $c(n,m)$, $c(n+1,m)$ and  $\pi-\alpha$ for
$c(n,m)$, $c(n,m+1)$ the governing equation (\ref{ortho-sg})
becomes
\begin{equation} \label{skew-sg}
R^2(r_1+r_2+r_3+r_4)-(r_2r_3r_4+r_1r_3r_4+r_1r_2r_4+r_1r_2r_3)+2R\cos \alpha
(r_1r_3-r_2r_4)=0.
\end{equation}
It is easy to see that (\ref{skew-sg}) has the same solution
(\ref{erf-r}) and it therefore defines a square grid circle
pattern, which is a discrete $\rm Erf$. A hexagonal analog of $\rm
Erf$ is not known.

\subsection{Circle patterns with quasi-regular combinatorics.}

 One can deregularize the
prescribed combinatorics by a projection of ${\mathbb Z}^n$ into a
plane as follows (see \cite{Se}). Consider ${{\mathbb
Z}^n_+}\subset {\mathbb R}^n$. For each coordinate vector ${\bf
e}_i=(e^1_i,...,e^n_i)$, where $e^j_i=\delta ^j_i$ define a unit
vector ${\bf \xi}_i$ in ${\mathbb C}={\mathbb R}^2$ so that for
any pair of indices $i,j$, vectors ${\bf \xi}_i,{\bf \xi}_j$ form
a basis in ${\mathbb R}^2$. Let $\Omega\in{\mathbb R}^n$ be some
2-dimensional simply connected cell complex with vertices in
${\mathbb Z}^n_+$. Choose some ${x_0}\in \Omega$.
Define
the map $P:\Omega \to {\mathbb C}$ by the following conditions:
\begin{itemize}
\item$P(x_0)=P_0$,
\item if  $x,y$ are vertices of $\Omega$ and $y=x+{\bf e}_i$ then
$P(y)=P(x)+{\bf \xi}_i$.
\end{itemize}

\noindent It is easy to see that $P$ is correctly defined and
unique. \\
We call $\Omega$ a {\it projectable} cell complex if its image
$\omega =P(\Omega)$ is embedded, i.e. intersections of images of
different cells of $\Omega$ do not have inner parts. Using
projectable cell complexes one can obtain combinatorics of Penrose
tilings.

\medskip

\noindent It is natural to define ``discrete conformal map on
$\omega$'' as a discrete complex immersion function $z$ on
vertices of $\omega$ preserving the cross-ratios of the
$\omega$-cells. The argument of $z$ can be labeled by the vertices
$x$ of $\Omega$. Hence for any cell of $\Omega$, constructed on
${\bf e}_k, {\bf e}_j$, the function $z$ satisfies the following
equation for the cross-ratios:
\begin{equation} \label{g-q}
q(z_x,z_{x+{\bf e_k}},z_{x+{\bf e_k}+{\bf e_j}},z_{x+{\bf
e_j}})=e^{-2i\alpha_{k,j} },
\end{equation}
where $\alpha_{k,j} $ is the angle between ${\bf \xi}_k$ and ${\bf
\xi}_j$, taken positively if $({\bf \xi}_k, {\bf \xi}_j)$ has
positive orientation and taken negatively otherwise. Now suppose
that $z$ is a solution to (\ref{g-q}) defined on the whole
${\mathbb Z}^n_+$. We can define a discrete $z^c:\omega \to
{\mathbb C }$ for projectable $\Omega$ as a solution to
(\ref{g-q}),(\ref{g-c}) restricted to $\Omega$. Initial conditions
for this solution are of the form (\ref{initial}) so that the
restrictions of $z$ to each two-dimensional coordinate lattice is
an immersion defining a circle pattern with prescribed
intersection angles.
This
definition naturally generalizes the definition of discrete
hexagonal and square grid $z^c$ considered above.

We finish this section with the natural conjecture formulated in
\cite{A1}.
\medskip

\noindent {\bf Conjecture.} {\it The discrete $z^c:\omega \to
{\mathbb C}$ is an immersion.}
\medskip

\noindent First step in proving this claim is to show that
equation (\ref{g-q}) is compatible with the constraint
\begin{equation} \label{g-c}
c f_x=\sum_{s=1}^n2x_s \frac{(f_{x+{\bf e}_s}-f_x)(f_x-f_{x-{\bf
e}_s})}{f_{x+{\bf e}_s}-f_{x-{\bf e}_s}}
\end{equation}
For $n=3$ this fact is proven in \cite{BH}.

\section{Acknowledgements}

This research was supported by the DFG Research Center
"Mathematics for key technologies" (FZT 86) in Berlin
 and the EPSRC grant No Gr/N30941.


\begin{thebibliography}{99}

\addcontentsline{toc}{section}{References}

\bibitem{AB}
S.I.Agafonov, A.I.Bobenko,  Discrete $Z^\gamma$ and Painlev\'e
equations, {\it Internat. Math. Res.
 Notices,}  {\bf 4} (2000), 165-193

\bibitem{A} S.I.Agafonov, Embedded circle patterns with the
combinatorics of the square grid  and discrete Painlev\'e
equations, {\it Discrete Comput. Geom.,} {\bf 29:2} (2003),
305-319

\bibitem{A1} S.I.Agafonov, Discrete Riccati equation, hypergeometric functions and circle
patterns of Schramm type, Preprint arXiv:math.CV/0211414


\bibitem{BDS}
A.F.Beardon, T.Dubejko, K.Stephenson,  Spiral hexagonal circle
packings in the plane,   {\it Geom. Dedicata,} {\bf 49} (1994),
39-70

\bibitem{BS}
A.F.Beardon, K.Stephenson  The uniformization theorem for circle packings,
{\it Indiana Univ. Math. J.,} {\bf 39:4} (1990), 1383-1425

\bibitem{BHConf} A.I.Bobenko, T.Hoffmann,
 Conformally symmetric circle packings. A generalization of Doyle spirals,
{\it Experimental Math.,}  {\bf 10:1} (2001), 141-150

\bibitem{BHS} A.I.Bobenko, T.Hoffmann, Yu.B.Suris,
 Hexagonal circle patterns and integrable systems:
Patterns with the multi-ratio property and Lax equations on the
regular triangular lattice, {\it Internat. Math. Res.
 Notices,}  {\bf 3} (2002), 111-164.

\bibitem{BH} A.I.Bobenko, T.Hoffmann,
 Hexagonal circle patterns and integrable systems. Patterns with constant angles, to appear in {\it Duke Math.
 J.,}
 Preprint arXiv:math.CV/0109018


 \bibitem{BPD}
A.~Bobenko, U.~Pinkall,   Discretization of surfaces
and integrable
systems,
 In: {\it Discrete Integrable Geometry and Physics;} Eds.
A.I.Bobenko and
R.Seiler,  pp. 3-58,   Oxford University Press, 1999.

\bibitem{Boole} G.~Boole, Calculus of finite differences, Chelsea
Publishing Company, NY 1872.

\bibitem{Doy}
K.~Callahan, B.~Rodin,   Circle packing immersions
form regularly
exhaustible surfaces, {\it Complex Variables,}
{\bf 21} (1993), 171-177.

\bibitem{DS} T.Dubejko, K.Stephenson,  Circle packings: Experiments in discrete
 analytic function theory, {\it Experimental Math.} {\bf 4:4}, (1995), 307-348

\bibitem{HS}
Z.-X.~He, O.~Schramm,   The $C^{\infty}$ convergence
of hexagonal disc
packings to Riemann map, {\it Acta. Math.}
{\bf 180} (1998),
219-245.

\bibitem{H} Z.-X.~He,   Rigidity of infinite disk patterns, {\it Annals
of
Mathematics} {\bf 149} (1999), 1-33.

\bibitem{TH} T.~Hoffmann,  Discrete CMC surfaces and
discrete
holomorphic maps,
In: {\it Discrete Integrable Geometry and Physics,}
Eds.: A.I.Bobenko and
R.Seiler, pp. 97-112,  Oxford University Press, 1999.


\bibitem{MR} A.Marden, B.Rodin,  On Thurston's formulation and proof of
 Andreev's theorem, pp. 103-115, Lecture Notes in Math., Vol. 1435, 1990.

\bibitem{NC} F.~Nijhoff, H.~Capel,  The discrete Korteveg de
Vries
Equation,
{\it Acta Appl. Math.} {\bf 39} (1995), 133-158.

\bibitem{N} F.W.~Nijhoff,
Discrete Painlev\'e equations and symmetry reduction on the lattice.
In: {\it Discrete Integrable Geometry and Physics;} Eds.
A.I.Bobenko and
R.Seiler,  pp. 209-234,   Oxford University Press, 1999.

\bibitem{R} B.~Rodin, Schwarz's lemma for circle packings, {\it Invent. Math.,} {\bf 89}
 (1987), 271-289.

\bibitem{RS}
B.~Rodin, D.~Sullivan,  The convergence of circle packings to
Riemann mapping, {\it J. Diff. Geometry} {\bf 26}
 (1987), 349-360.


\bibitem{Schramm} O.~Schramm,
Circle patterns with the
combinatorics of the square grid, {\it Duke Math. J.}
{\bf 86} (1997),
347-389.


\bibitem{Se} M.~Senechal, {\it Quasicrystals and geometry.} Cambridge University Press, Cambridge, 1995.


\bibitem{T} W.P.~Thurston,   The finite Riemann mapping
theorem,
Invited talk, International Symposium on the occasion of
the proof of the
Bieberbach Conjecture, Purdue University (1985).



\end{thebibliography}
\end{document}